\documentclass[a4paper,11pt]{amsart}
\usepackage{fullpage}
\usepackage{amsfonts}
\usepackage{amssymb}
\usepackage{amssymb}
\usepackage{amsmath,amsthm}
\usepackage[colorlinks=true,
            linkcolor=red,
            urlcolor=blue,
            citecolor=blue]{hyperref}
\usepackage{amsthm}
\usepackage[english]{babel}
\usepackage{mathrsfs}
\usepackage[mathscr]{eucal}
\usepackage{cite}
\usepackage{upgreek}

\numberwithin{equation}{section}

\newtheorem{main}{Theorem}

\newtheorem{mcor}[main]{Corollary}

\newtheorem{theorem}{Theorem}[section]
\newtheorem{lem}[theorem]{Lemma}
\newtheorem{prop}[theorem]{Proposition}

\newtheorem{cor}[theorem]{Corollary}

\theoremstyle{definition}

\newtheorem{assumption}[theorem]{Assumption}

\newtheorem{Remark}[theorem]{Remark}

\newtheorem{examples}[theorem]{Example}
\newtheorem{claim}[theorem]{Claim}
\newtheorem*{claim*}{Claim}
{ \theoremstyle{remark} }

\DeclareMathAlphabet{\pazocal}{OMS}{zplm}{m}{n}

\newcommand{\paM}{\pazocal{M}}

\def\ra{\rightarrow}

\def\La{\Lambda}
\def\Om{\Omega}
\def\Th{\Theta}

\def\G{\Gamma}

\theoremstyle{remark}

\newcommand{\set}[1]{ \{ #1 \} }

\begin{document}

 \title[Free Product Rigidity for Group Algebras]{Amalgamated Free Product Rigidity for \\ Group von Neumann Algebras}

\author[I. Chifan]{Ionu\c{t} Chifan}
\address{Department of Mathematics, The University of Iowa, 14 MacLean Hall, Iowa City, IA  
52242, USA}
\email{ionut-chifan@uiowa.edu}
\thanks{ I.C. was partially supported by NSF grants DMS \# 1600688 and DMS \# 1301370. Part of this work was done during a Flex Load Semester awarded by the CLAS at University of Iowa.}

\author[A. Ioana]{Adrian Ioana}
\address{Department of Mathematics, University of California San Diego, 9500 Gilman Drive, La Jolla, CA 92093, USA, and IMAR, Bucharest, Romania}
\email{aioana@ucsd.edu}
\thanks{A.I. was partially supported by  NSF Career Grant DMS \#1253402 and a Sloan Foundation Fellowship.}

\maketitle







\begin{abstract} 
We provide a fairly large family of amalgamated free product groups $\Gamma=\Gamma_1\ast_{\Sigma}\Gamma_2$ \\ whose amalgam structure can be completely recognized from their von Neumann algebras. Specifically, assume that $\Gamma_i$ is a product of two icc non-amenable bi-exact groups, and $\Sigma$ is icc amenable with trivial one-sided commensurator in $\Gamma_i$, for every $i=1,2$. Then $\Gamma$ satisfies the following rigidity property: \emph{any} group $\La$  such that $L(\La)$ is isomorphic  to $L(\G)$ admits an amalgamated free product decomposition $\La=\La_1\ast_\Delta \La_2$ such that the inclusions $L(\Delta)\subseteq L(\La_i)$ and $L(\Sigma)\subseteq L(\G_i)$ are isomorphic, for every $i=1,2$. This result significantly strengthens some of the previous Bass-Serre rigidity results for von Neumann algebras.
As a corollary, we obtain the first examples of amalgamated free product groups which are W$^*$-superrigid.
\end{abstract}

\maketitle


\section{Introduction}

\noindent In  \cite{MvN36, MvN43}, Murray and von Neumann found a natural way to associate a von Neumann algebra, denoted by $L(\G)$, to every countable discrete group $\G$. More precisely, $L(\G)$ is defined as the weak operator closure of the complex group algebra $\mathbb C\Gamma$ acting by left convolution on the Hilbert space $\ell^2(\G)$.
The classification of group von Neumann algebras has since been a central theme in operator algebras driven by the following question: what aspects of the group $\Gamma$ are remembered by $L(\Gamma)$?  This question is the most interesting when $\Gamma$ is {\it icc} (i.e., the conjugacy class of every non-trivial element of $\Gamma$ is infinite), which corresponds to $L(\Gamma)$ being a II$_1$ factor.

 Von Neumann algebras tend to forget a lot of information about the groups they are constructed from.
This is best illustrated by Connes'  
 theorem asserting that all II$_1$ factors arising from icc amenable groups are isomorphic to the hyperfinite II$_1$ factor \cite{Co76}. Consequently, amenable groups manifest a striking lack of rigidity: any algebraic property of the group (e.g., being torsion free or finitely generated) is completely lost in the passage to von Neumann algebras.

In sharp contrast, in the non-amenable case, Popa's deformation/rigidity theory has led to the discovery of several instances when various properties of a group $\Gamma$ can be recovered from $L(\Gamma)$. We only highlight three developments in this direction here, and refer the reader to the surveys \cite{Po06a,Va10,Io12} for more information.  Thus, it was shown in \cite{Po03,Po04} that within a large of icc groups (containing the wreath product $\mathbb Z/2\mathbb Z\wr\G$, for any infinite property (T) group $\Gamma$) isomorphism of the associated II$_1$ factors implies isomorphism of the groups. 
A few years later,
 the first examples of  groups, called \emph{$W^*$-superrigid} groups, that can be entirely reconstructed from their von Neumann algebras were discovered in \cite{IPV10} (see\cite{BV13,Be14} for the only other know examples).
Specifically, a group $\Gamma$ is called W$^*$-superrigid if whenever $L(\Lambda)$ is isomorphic to $L(\Gamma)$, $\Lambda$ must be isomorphic to $\Gamma$. 
Most recently, the following \emph{product rigidity} phenomenon was found in \cite{CdSS15}: if $\G_1$ and $\G_2$ are icc hyperbolic groups, then any group $\La$ such that $L(\La)$ is isomorphic to $L(\G_1\times\G_2)$ admits a decomposition $\La=\La_1\times \La_2$ such that  $L(\La_i)$ is isomorphic to $L(\G_i)$, up to amplifications, for every $i\in\{1,2\}$.  In other words,  the von Neumann algebra $L(\Gamma)$ completely remembers the product structure of the underlying group $\Gamma$.

Motivated by these advances, it seems natural to investigate instances when other constructions in group theory 
can be recognized from the von Neumann algebraic structure. We make progress on this general problem here by providing a class of amalgamated free product (abbreviated AFP) groups $\Gamma$ whose von Neumann algebra $L(\Gamma)$ entirely remembers the amalgam structure of $\Gamma$.

Before stating our main result, let us recall  the definition of bi-exact groups \cite[Definition 15.1.2]{BO08}. 
A countable group $\Gamma$ is said to be {\it bi-exact} (or to belong to Ozawa's class $\mathcal S$ \cite{Oz04}) if it is exact and admits a map $\mu:\Gamma\rightarrow\text{Prob}(\Gamma)$ satisfying $\lim_{x\rightarrow\infty}\|\mu(gxh)-g\cdot\mu(x)\|=0$, for all $g,h\in\Gamma$. The class of bi-exact groups includes all hyperbolic groups \cite{Oz03}, the wreath product $A\wr\Gamma$ of any amenable group $A$ with a bi-exact group $\Gamma$ \cite{Oz04}, the  group $\mathbb Z^2\rtimes SL_2(\mathbb Z)$ \cite{Oz08}, and is closed under free products.

\begin{main}\label{A} 
Let $\G=\Gamma_1\ast_{\Sigma}\Gamma_2$ be an amalgamated free product group satisfying the following:

\begin{enumerate} 
\item $\Sigma$ is an icc amenable group and $[\Sigma:\Sigma\cap g\Sigma g^{-1}]=\infty$, for every $g\in\Gamma_i\setminus\Sigma$ and $i\in\{1,2\}$.
\item $\Gamma_i=\Gamma_i^1\times\Gamma_i^2$, where $\Gamma_i^j$ is an icc, non-amenable, bi-exact group, for every $i,j\in\{1,2\}$. 
\end{enumerate}

Denote $M=L(\Gamma)$ and let $\Lambda$ be an arbitrary group such that $M=L(\Lambda)$.
 
 Then there exist a decomposition $\La=\La_1\ast_{\Delta}\La_2$ and a unitary $u\in M$ such that 
$$
 u L(\La_1)u^*=L(\Gamma_1),\;\;
uL(\La_2)u^*=L(\Gamma_2),\;\;
\text{and}\;\; uL(\Delta)u^*=L(\Sigma). 
$$
\noindent

\end{main}

Before placing Theorem \ref{A} into context, let us present some classes of groups to which it applies.

\begin{examples}\label{ex} Let $\Sigma_0<\Gamma_0$ be an inclusion of groups satisfying the following condition: $(\star)$
$\Gamma_0$ is icc, non-amenable, bi-exact, $\Sigma_0$ is icc, amenable, and $[\Sigma_0:\Sigma_0\cap g\Sigma_0 g^{-1}]=\infty$, for all $g\in\Gamma_0\setminus\Sigma_0$. 
Having such an inclusion $\Sigma_0<\Gamma_0$, the hypothesis of Theorem \ref{A} is satisfied for $\Gamma_1=\Gamma_2=\Gamma_0\times\Gamma_0$, and $\Sigma<\Gamma_1\cap\Gamma_2$ equal to either $\Sigma_0\times\Sigma_0$ or $\Sigma=\{(g,g)|g\in\Sigma_0\}$. 

On the other hand, $(\star)$ is verified by the following group inclusions:
\begin{enumerate}
\item[{(a)}] $A<A*B$, where $A$ is any  icc amenable group, and $B$ is any non-trivial bi-exact group.
\item[{(b)}] $A\wr C<A\wr D$, where $A$ is any non-trivial amenable group and $C$ is any infinite maximal amenable subgroup of any icc hyperbolic group $D$ (see Section \ref{exb} for a proof of this assertion). For instance, take $C=\mathbb Z$ and $D=C*\mathbb F_{n-1}=\mathbb F_n$, for some $2\leqslant n\leqslant+\infty$.

  \item[{(c)}] $\mathbb Z^2\rtimes\langle M\rangle<\mathbb Z^2\rtimes\mathbb F_n$, where we view $\mathbb F_n$ as a subgroup of $SL_2(\mathbb Z)$, for some $2\leqslant n\leqslant+\infty$, and $M\in\mathbb F_n$ is any matrix such that $|\text{Tr}(M)|>2$ and $M\not=M_0^{\ell}$, for every $M_0\in\mathbb F_n$ and $\ell\geqslant 2$.
    For instance, define $\mathbb F_2=\langle M_1,M_2\rangle<SL_2(\mathbb Z)$ as the group generated by 
    $M_1=\begin{pmatrix}1&2\\0&1 \end{pmatrix}, M_2=\begin{pmatrix}1&0\\2&1\end{pmatrix}$  and let $M=M_1M_2=\begin{pmatrix}5&2\\2&1 \end{pmatrix}$.
\end{enumerate}
\end{examples}

The fact that a large class of AFP groups satisfy Theorem \ref{A} should not be surprising since II$_1$ factors of such groups (and, more generally, AFP II$_1$ factors $M=M_1\ast_{B}M_2$) have been shown to be   
extremely rigid  (see, e.g., \cite{Oz04, IPP05, Pe06, Po06b, CH08, Ki09, PV09, HPV09, Io12, Va13, HU15}). 
In particular,  Bass-Serre-type rigidity results for AFP II$_1$ factors have been discovered in \cite{IPP05}. 
For instance, assume  that $\Gamma=\Gamma_1\ast_{\Sigma}\Gamma_2$ and $\La=\La_1*_{\Delta}\La_2$ are AFP groups, where  $\Gamma_1,\Gamma_2,\Lambda_1,\Lambda_2$ are icc property (T) groups. 
It follows from \cite{IPP05} that if $\theta:L(\G)\ra L(\La)$ is any $*$-isomorphism, 
then $\theta(L(\G_i))$ is unitarily conjugate to $L(\La_i)$, up to a permutation of indices.
Later on, the same conclusion was shown to hold assuming that $\Gamma_1,\Gamma_2,\Lambda_1,\Lambda_2$ are products of icc non-amenable groups and $\Sigma,\Delta$ are amenable \cite{CH08} (see also \cite{Oz04,Pe06} in the case $\Sigma$ and $\Delta$ are trivial).

Theorem \ref{A} strengthens such Bass-Serre rigidity results in the case $\Sigma$ is icc, by removing all assumptions on the group $\La$.
This is strongest type of rigidity that one can expect for II$_1$ factors of general AFP groups.
To make this precise, note that if $\Gamma=\Gamma_1*_{\Sigma}\Gamma_2$, then $L(\G)$ is determined up to isomorphism by the isomorphism classes of the inclusions $L(\Sigma)\subseteq L(\G_1)$ and $L(\Sigma)\subseteq L(\G_2)$. 
Conversely, Theorem \ref{A} asserts that, under certain assumptions, these isomorphism classes can be reconstructed from the isomorphism class of $L(\G)$.

\begin{Remark}
 We do not know whether Theorem \ref{A} holds for {\it plain} free product groups, i.e. when $\Sigma=\{e\}$. Note in this respect that by a result in \cite{DR01} isomorphism of the free group factors would imply that $L(\G_1*\G_2)\cong L(\G_1*\G_2*\mathbb F_{\infty})$, for any icc groups $\G_1$ and $\G_2$. However, even in the case $\Sigma=\{e\}$, we can still sometimes deduce that any group $\La$ with $L(\G)=L(\La)$ must contain subgroups  $\La_1, \La_2$ such that $L(\G_i)$ is unitarily conjugate to $L(\La_i)$, for every $i\in\{1,2\}$. 
 Indeed, in the context of Theorem \ref{A}, this holds if $\Gamma_1$ has property (T) (or Haagerup's property) while $\Gamma_2$ does not
 (see Corollary \ref{BS}).
\end{Remark}

\begin{Remark}
The groups covered by Theorem \ref{A} are typically not W$^*$-superrigid. 
Indeed, the  groups $(A*B\times A*B)*_{A\times A}(A*B\times A*B)$, where $A$ and $B$ are any icc amenable groups,
satisfy the hypothesis of Theorem \ref{A} (see Example \ref{ex}(a))  but produce isomorphic II$_1$ factors by \cite{Co76}.
\end{Remark}

Nevertheless, by combining Theorem \ref{A} with results from \cite{IPV10,CdSS15} we prove the following:

\begin{mcor}\label{B}
Let $\Gamma_0$ be an icc, non-amenable, bi-exact group, and $\Sigma_0<\Gamma_0$ be an icc, amenable subgroup such that the following two conditions hold:

\begin{enumerate}
\item $[\Sigma_0:\Sigma_0\cap g\Sigma_0 g^{-1}]=\infty$, for every $g\in\Gamma_0\setminus\Sigma_0$. 
\item the centralizer in $\Gamma_0$ of any finite index subgroup of $\Sigma_0\cap g\Sigma_0 g^{-1}$ is trivial, for any $g\in\Gamma_0$.
\end{enumerate}
Define $\Sigma:=\{(g,g)|g\in\Sigma_0\}<\Gamma_0\times\Gamma_0$ and $\G:=(\Gamma_0\times\Gamma_0)\ast_{\Sigma}(\Gamma_0\times\Gamma_0)$.

If $\La$ is any countable group and $\theta:L(\Gamma)\rightarrow L(\La)$ is any $*$-isomorphism, then there exist a group isomorphism $\delta:\G\ra\La$, a unitary $u\in L(\La)$, and a character $\eta:\G\ra\mathbb T$ such that $$\theta(u_g)=\eta(g) u v_{\delta(g)} u^*,\;\;\;\text{for every $g\in \G$.}$$

Here, $\{u_g\}_{g\in\G}$ and $\{v_h\}_{h\in\La}$ denote the canonical unitaries generating $L(\G)$ and $L(\La)$.
\end{mcor}

Corollary \ref{B} provides a new class of W$^*$-superrigid groups. Note that unlike all of the known classes of  W$^*$-superrigid groups, our examples are not generalized wreath product groups nor special subgroups of such groups, as in \cite{IPV10,BV13,Be14}. As such, Corollary \ref{B} gives first the examples of AFP groups which are W$^*$-superrigid.

\begin{examples}
To give some concrete examples of inclusions $\Sigma_0<\Gamma_0$ to which Corollary \ref{B} applies, let $A$ be any non-trivial amenable group and $C$ be any infinite maximal amenable subgroup of any icc hyperbolic group $D$. Put $\Sigma_0:=A\wr C$ and $\Gamma_0:=A\wr D$.
Then condition (1) from Corollary \ref{B} is satisfied by Example \ref{ex}(b). If $g\in\Gamma_0$, then $\Sigma_0\cap g\Sigma_0g^{-1}\supseteq A^{(C)}=\oplus_{c\in C}A$. Thus, any finite index subgroup of $\Sigma_0\cap g\Sigma_0g^{-1}$ contains $A_0^{(C)}$, for some finite index subgroup $A_0<A$. Since $A$ is icc, the centralizer of $A_0^{(C)}$ in $\Gamma_0$ is trivial, which proves that condition (2) from Corollary \ref{B} is also satisfied.

In particular, Corollary \ref{B} applies in the case $\Sigma_0=S_{\infty}\wr\mathbb Z$ and $\Gamma_0=S_{\infty}\wr\mathbb F_n$,  for some $n\geqslant 2$, where $S_{\infty}$ denotes the group of finite permutations of $\mathbb N$.
\end{examples}

We end by noticing that the groups $\Gamma$ from Corollary \ref{B} are also {\it C$^*$-superrigid}, in the sense that they can be completely recovered from their reduced C$^*$-algebras. Recall in this respect that the reduced C$^*$-algebra C$^*_{\text{r}}(\Gamma)$ of $\Gamma$ is defined as the operator norm closure of the linear span of the left regular representation $\{u_g\}_{g\in\Gamma}\subseteq\mathcal U(\ell^2(\Gamma))$. A countable group $\Gamma$ is then called C$^*$-superrigid if whenever $C^*_r(\Lambda)$ is isomorphic to $C^*_r(\Gamma)$, $\Lambda$ must be isomorphic to $\Gamma$.

\begin{mcor}\label{C}
Let $\G$ be any AFP group as in Corollary \ref{B}.

If $\La$ is any countable group and $\theta:\text{C}^*_{\text r}(\Gamma)\rightarrow \text{C}^*_{\text{r}}(\La)$ is any $*$-isomorphism, then there exist a group isomorphism $\delta:\G\ra\La$, a unitary $u\in L(\La)$, and a character $\eta:\G\ra\mathbb T$ such that $$\theta(u_g)=\eta(g) u v_{\delta(g)} u^*,\;\;\;\text{for every $g\in \G$.}$$
\end{mcor}
The first examples of non-abelian torsion-free C$^*$-superrigid groups were recently found in 
 \cite{KRTW}. The groups considered in \cite{KRTW} are virtually abelian, hence amenable. In contrast, Corollary \ref{C} provides the first examples of non-amenable groups that are C$^*$-superrigid.

\subsection*{Comments on the proof of Theorem \ref{A}} 
We end the introduction with some brief and informal comments on the proof of Theorem \ref{A}.
Let $\Gamma=\Gamma_1*_{\Sigma}\Gamma_2$ be as in the hypothesis of Theorem \ref{A}, where $\Gamma_i=\Gamma_i^1\times\Gamma_i^2$ is a product of icc, non-amenable, bi-exact groups, for every $i\in\{1,2\}$. Our goal is to investigate all possible group von Neumann algebra decompositions of $M=L(\G)$. 
To this end, let $\La$ be a countable group such that $M=L(\La)$, and consider the $*$-homomorphism $\triangle:M\ra M\bar{\otimes}M$ given by $\triangle(v_h)=v_h\otimes v_h$, for all $h\in\La$ \cite{PV09}. 

In the first part of the proof, we use Ozawa's work on subalgebras with non-amenable commutant inside von Neumann algebras of (relatively) bi-exact groups \cite{Oz03,Oz04, BO08} to conclude that \begin{equation}\label{triangle}\text{$\triangle(L(\Gamma_1^1))\prec M\bar{\otimes}L(\Gamma_m^n)$, for some $m,n\in\{1,2\}$}.\end{equation} Here, $P\prec Q$ denotes the fact that a corner of $P$ embeds into a corner of $Q$ inside the ambient algebra, in the sense of Popa \cite{Po03}.

The second part of the proof consists of combining \eqref{triangle} with an ultrapower technique from \cite{Io11} to deduce the existence of a subgroup $\Om<\La$ such that  \begin{equation}\label{triangle2}\text{$L(\G_1^1)\prec L(\Om)$ and the centralizer of $\Om$ in $\La$ is non-amenable}.\end{equation} 
For these first two parts, see Theorem \ref{comm1}.

Further, by using \eqref{triangle2}  and building on techniques  from \cite{IPP05,CdSS15}  we find a subgroup $\Th<\La$  such that  
$L(\Th)z$ and $L(\Th)(1-z)$ are unitarily conjugate to corners of $L(\G_1)$ and $L(\G_2)$,
respectively, for some non-zero central projection $z\in L(\Th)$ (see Theorems \ref{peripheralidentification1} and \ref{peripheralidentification2}). More precisely, $\Th$ is defined as the one-sided commensurator of $\Om$ in $\La$ \cite{FGS10}.

The final part of the proof is the subject of Section \ref{4}. 
Thus, we first use that $\Sigma$ is icc to derive that $z=1$ (see Proposition \ref{noniso}). As a consequence of this and the analogous analysis for $\Gamma_2^1$ instead of $\Gamma_1^1$, we obtain subgroups $\Th_1,\Th_2<\La$ such that $L(\Th_i)$ is unitarily conjugate to $L(\G_i)$. With some additional work, we finally prove that $\La=\Th_1*_{\Th_1\cap h\Th_2h^{-1}}(h\Th_2h^{-1})$, for some $h\in\La$, and that the conclusion holds for $\La_1=\Th_1$ and $\La_2=h\Th_2h^{-1}$. 


\section{Preliminaries}

\noindent We begin this section by reviewing several concepts in von Neumann algebras and group theory. We then record several technical ingredients for the proofs of our main results.

\subsection{Terminology} 
All von Neumann algebras $M$ considered in this article are tracial, i.e., they are endowed with a unital, faithful, normal linear functional $\uptau:M\ra \mathbb C$  satisfying $\uptau(xy)=\uptau(yx)$, for all $x,y\in M$. Given $x\in M$, we denote by $\|x\|$ its operator norm, and by $\|x\|_2=\uptau(x^*x)^{1/2}$ its so-called $2$-norm.
For a tracial von Neumann algebra $M$, we denote by $\mathcal U(M)$ its unitary group, by $\mathcal Z(M)$ its center,  by $\mathcal P(M)$ the set of its projections, and by $(M)_1=\{x\in M\;|\;\|x\|\leq 1\}$ its unit ball with respect to the operator norm. 
For a set $S\subseteq M$, we denote by W$^*(S)$ the smallest von Neumann subalgebra of $M$ which contains $S$.
If $S$ is closed under adjoint, then W$^*(S)$ is equal to the bicommutant $S''$ of $S$ .

All inclusions $P\subseteq M$ of von Neumann algebras are assumed unital, unless otherwise specified. 
Given von Neumann subalgebras $P, Q\subseteq M$,  we denote by $E_{P}:M\rightarrow P$ the conditional expectation onto $P$, 
by $P'\cap M=\{x\in M\;|\;xy=yx,\;\text{for all $y\in P$}\}$ the relative commutant of $P$ in $ M$,
and by $P\vee Q=W^*(P\cup Q)$ the von Neumann algebra generated by $P$ and $Q$.

All groups considered in this article are countable and discrete.
For a group $\G$, we denote by $\{ u_g\}_{g\in\G} \subseteq \mathcal U(\ell^2(\G))$ its left regular representation given by $u_g(\delta_h ) = \delta_{gh}$, where $\delta_h$ is the Dirac mass at $h$. The weak operator closure of the linear span of $\{ u_g \}_{g\in G}$ in $\mathbb B(\ell^2(\G))$ is called the group von Neumann algebra of $\Gamma$ and is denoted by $L(\G)$. $L(\G)$ is a II$_1$ factor precisely when $\G$ has infinite non-trivial conjugacy classes (icc) \cite{MvN43}.

 Let $\G$ be a group. Given subsets $K,F\subseteq \G$, we put $KF=\{ kf | k\in K,f\in F\}$. Given a subgroup $\Sigma<\G$, we denote by $C_{\Sigma}(K)=\{ g \in \Sigma | g k=k g, \text{ for all }k\in K\}$  the centralizer of $K$ in $\Sigma$. For a positive integer $n$, we denote by $\overline{1,n}$ the set $\{1,2,...,n\}$.


\subsection{Popa's intertwining technique} In the early 2000s, S. Popa introduced in \cite [Theorem 2.1 and Corollary 2.3]{Po03} the following powerful criterion for the existence of intertwiners between arbitrary subalgebras of tracial von Neumann algebras.

\begin {theorem}[\!\!\cite{Po03}]\label{corner} Let $(M,\uptau)$ be a separable tracial von Neumann algebra and let $ P,  Q\subseteq  M$ be (not necessarily unital) von Neumann subalgebras. 
Then the following are equivalent:
\begin{enumerate}
\item There exist $ p\in  \mathcal P( P), q\in  \mathcal P( Q)$, a $\ast$-homomorphism $\theta:p  P p\rightarrow q Q q$  and a non-zero partial isometry $v\in q M p$ such that $\theta(x)v=vx$, for all $x\in p P p$.
\item For any group $\mathcal U\subset \mathscr U( P)$ such that $\mathcal U''=  P$ there is no sequence $(u_n)_n\subset \mathcal U$ satisfying $\|E_{  Q}(xu_ny)\|_2\rightarrow 0$, for all $x,y\in   M$.
\end{enumerate}
\end{theorem} 

 If one of the two equivalent conditions from Theorem \ref{corner} holds we say that \emph{ a corner of $ P$ embeds into $ Q$ inside $ M$}, and write $ P\prec_{ M} Q$. If we moreover have that $ P p'\prec_{ M} Q$, for any nonzero projection  $p'\in  P'\cap 1_ P  M 1_ P$, then we write $ P\prec_{ M}^{s} Q$. 

Next, we record two basic intertwining results that will be used later on.


\begin{lem}\label{intertwiningroups} Let $\Gamma_1,\Gamma_2<\Gamma$ be countable groups such that $L(\Gamma_1)\prec_{L(\Gamma)} L(\Gamma_2)$.

 Then  there exists $g\in \G$ such that $[\G_1:\G_1\cap g\G_2g^{-1}]<\infty$. 
\end{lem}

{\it Proof.}
 Denote by $e$ the orthogonal projection from $\ell^2(\Gamma)$ onto $\ell^2(\Gamma_2)$. Consider Jones' basic construction $\langle L(\Gamma), e\rangle\subseteq\mathbb B(\ell^2(\Gamma))$ of the inclusion $L(\G_2)\subseteq L(\G)$ endowed with the usual semi-finite trace $Tr$ and 2-norm $\|x\|_{2,Tr}=Tr(x^*x)^{1/2}$ (see e.g. \cite{Jo81,PP86}).  Denote by $\mathcal H:=L^2(\langle L(\Gamma),e\rangle)$ the Hilbert space obtained by completing $\langle L(\Gamma),e\rangle$ with respect to $\|.\|_{2,Tr}$. Let $\pi:\Gamma\rightarrow\mathcal U(\mathcal H)$ denote the unitary representation given by $\pi(g)(\xi)=u_g\xi u_g^*$.  

\begin{claim}\label{invariant}
If $\xi\in\mathcal H$ is a $\pi(\Gamma_1)$-invariant vector, then $\xi$ belongs to the $\|.\|_{2,Tr}$-closure of the linear span of $\{u_geu_h\;|\;g,h\in\Gamma\;\text{with}\;[\Gamma_1:\Gamma_1\cap g\Gamma_2g^{-1}]<\infty\}$.
\end{claim}

{\it Proof of Claim \ref{invariant}.}
Let $\mathcal H_k\subseteq\mathcal H$ be the $\|.\|_{2,Tr}$-closure
 of the linear span of $\{u_{gk}eu_g^*|g\in\Gamma\}$. Then $\mathcal H_k$ is $\pi(\Gamma)$-invariant, for every $k\in\Gamma$, and if $S\subseteq\Gamma$ is such that $\Gamma=\sqcup_{k\in S}\{hkh^{-1}|h\in\Gamma_2\}$, then $\mathcal H=\bigoplus_{g\in S}\mathcal H_k$. 
 Thus, in order to prove the claim, we may assume that $\xi$ belongs to $\mathcal H_k$, for some fixed $k\in\Gamma$.
 Notice that
  $$\langle u_gu_keu_g^*,u_ke\rangle=\begin{cases}1,\;\;\text{if $g\in C_{\Gamma_2}(k)$}\\ 0,\;\;\text{if $g\notin C_{\Gamma_2}(k)$} \end{cases}$$  
 Thus, $\{u_{gk}eu_g^*\}_{g\in\Gamma/C_{\Gamma_2}(k)}$ is a $\pi(\Gamma)$-invariant orthonormal basis of $\mathcal H_k$. 
Let $\xi\in\mathcal H_k$ be a $\pi(\Gamma_1)$-invariant vector and write $\xi=\sum_{g\in\Gamma/C_{\Gamma_2}(k)}c_gu_{gk}eu_g^*$, for some scalars $c_g$. If $c_g\not=0$, for some $g\in\Gamma/C_{\Gamma_2}(k)$, then $\pi(\Gamma_1)(u_{gk}eu_g^*)$ must be finite, or equivalently $[\Gamma_1:\Gamma_1\cap gC_{\Gamma_2}(k)g^{-1}]<\infty$. 
 This implies that $[\Gamma_1:\Gamma_1\cap (gk)\Gamma_2(gk)^{-1}]<\infty$, which yields the claim.
 \hfill$\square$

We are now ready to derive  the lemma.
Since $L(\G_1)\prec_{L(\Gamma)}L(\G_2)$, \cite[Theorem 2.1]{Po03} implies that the $L(\Gamma)$-$L(\Gamma)$-bimodule $\mathcal H$ contains a non-zero $L(\G_1)$-central vector. Thus, $\mathcal H$ contains a non-zero $\pi(\Gamma_1)$-invariant vector, and the Claim \ref{invariant} implies the conclusion.
\hfill$\blacksquare$

\begin{lem}\label{intersect} Let $\G_1,\G_2<\G$ be countable groups, and $ Q\subseteq qL(\G)q$ be a von Neumann subalgebra. For every $i\in\overline{1,2}$, suppose that $p_i\in\mathcal P( Q'\cap qL(\G)q)$ and  $u_i\in\mathcal U(L(\G))$ satisfy $u_i Q p_iu_i^*\subseteq L(\G_i)$. Assume that $p_1p_2\not=0$ and let $p\in  Q'\cap qL(\G)q$ be a projection such that $pp_1p_2\not=0$.

\vskip 0.03in
\noindent
Then there exists $g\in\Gamma$ such that $ Q p\prec_{L(\G)} L(\Gamma_1\cap g\Gamma_2g^{-1})$.

Moreover, if $Q\subseteq L(\G_1)\cap vL(\G_2)v^*$, for some partial isometry $v\in L(\G)$ satisfying $vv^*=q$ and $v^*v\in L(\G_2)$,   then we can find $g\in\G$  such that  $ Q\prec_{L(\G)} L(\Gamma_1\cap g\Gamma_2g^{-1})$ and   $\uptau(vu_g^*)\not=0$.
\end{lem}

{\it Proof.}
Assume that $pp_1p_2\not=0$, and put $\delta:=\|pp_1p_2\|_2>0$. For a set $S\subseteq\Gamma$, we denote by $e_S$ the orthogonal projection from $\ell^2(\Gamma)$ onto the closed linear span of $\{u_g|g\in S\}$. Note that \begin{equation}\label{belong}upp_1p_2\in pp_1u_2^*(L(\Gamma_2))_1u_2\cap pu_1^*(L(\Gamma_1))_1u_1p_2,\;\;\;\text{for all $u\in\mathcal U(Q)$.}\end{equation}

Let $S_1\subseteq \Gamma$ be a finite set such that  $\|pp_1u_2^*-v_1\|_2\leqslant\delta/5$, where $v_1=e_{S_1}(pp_1u_2^*)$. By Kaplansky's density theorem, for any $i\in\overline{2,4}$, we can find  $S_i\subseteq\Gamma$ finite and $v_i\in (L(\Gamma))_1$ belonging to the linear span of $\{u_g|g\in S_i\}$ such that $\|v_2-u_2\|_2\leqslant \delta/5$, $\|v_3-pu_1^*\|_2\leqslant\delta/5$, and $\|v_4-u_1p_2\|_2\leqslant\delta/5$.
 By using these inequalities and \eqref{belong} we get that $\|upp_1p_2-e_{S_1\Gamma_2S_2\cap S_3\Gamma_2 S_4}(upp_1p_2)\|_2\leqslant 4\delta/5$, for all $u\in\mathcal U(Q)$. Since $\|upp_1p_2\|_2=\delta$, we get that 
$$
\|e_{S_1\Gamma_2S_2\cap S_3\Gamma_1S_4}(upp_1p_2)\|_2\geqslant\delta/5,\;\;\;\text{for every $u\in\mathcal U( Q)$.}
$$
It is easy to see that there is $S\subseteq\Gamma$ finite such that $S_1\Gamma_2S_2\cap S_3\Gamma_1S_4\subseteq \cup_{g,h,k\in S}h(\Gamma_1\cap g\Gamma_2g^{-1})k$.
Then the last inequality implies that $$\sum_{g,h,k\in S}\|E_{L(\Gamma_1\cap g\Gamma_2g^{-1})}(u_h^*upp_1p_2u_k^*)\|_2^2\geqslant \delta^2/25,\;\;\;\text{for every $u\in\mathcal U( Q)$.}$$
\noindent
By the proof of \cite[Theorem 4.3]{IPP05} (see also \cite[Remark 2.3]{DHI16}) this implies the conclusion. 

For the moreover assertion, put $\delta:=\|q\|_2$. Let $T_1,T_2\subseteq \G$ finite such that  $\|v-w_1\|_2\leqslant\delta/3$ and $\|v^*-w_2\|_2\leqslant\delta/3$, where 
$w_1=e_{T_1}(v)$ and $w_2\in (L(\G))_1$ satisfies $w_2=e_{T_2}(w_2)$. We may moreover assume that $\uptau(vu_g^*)\not=0$, for all $g\in T_1$. Then as above we find that $\|e_{\Gamma_1\cap T_1\G_2T_2}(u)\|_2\geqslant\delta/3$, for all $u\in\mathcal U(Q)$. Since $\G_1\cap T_1\G_2T_2\subseteq\cup_{g\in T_1,h\in T}(\G_1\cap g\G_2g^{-1})h$, for some finite set $T\subseteq\G$, we conclude that $Q\prec_{L(\G)}L(\G_1\cap g\G_2g^{-1})$, for some $g\in T_1$. This finishes the proof.
\hfill$\blacksquare$

\subsection{Commensurators and quasinormalizers} Let  $\Sigma< \G$ be an inclusion of countable groups, and $ P \subseteq  M$ be an inclusion of tracial von Neumann algebras.

 The \emph{commensurator} $\text{Comm}_\G(\Sigma)$ of $\Sigma$ in $\G$ is defined as the subgroup of all $g \in \G$ for which there exists a finite set $F\subseteq \G$ such that $\Sigma g\subseteq F\Sigma$ and $g\Sigma\subseteq \Sigma F$. 
Thus, $g\in \text{Comm}_\G(\Sigma)$ if and only if  $[\Sigma:\Sigma\cap g\Sigma g^{-1}]<\infty$ and $[g\Sigma g^{-1}:\Sigma\cap g\Sigma g^{-1}]<\infty$.
The {\it quasi-normalizer} $q\mathcal{N}_M( P)$ of $ P$ in $ M$ is defined as the $*$-algebra of all $x\in  M$ for which there exist $x_1,x_2,...,x_k\in  M$ such that  $ P x\subseteq \sum_i x_i  P$ and $x P \subseteq \sum_i  P x_i$ (see \cite[Definition 4.8]{Po99}).

In this paper, we will use the following one sided versions of these notions considered in \cite{FGS10}.
The \emph{one sided commensurator} $\text{Comm}^{(1)}_{\G}(\Sigma)$ is defined as the semigroup of all $g\in \G$ for which there exists a finite set $F\subseteq \G$ such that $\Sigma g\subseteq F\Sigma$. Thus, $g\in \text{Comm}^{(1)}_{\G}(\Sigma)$ if and only if $[\Sigma:\Sigma\cap g\Sigma g^{-1}]<\infty$.
 The {\it one sided quasi-normalizer} $q\mathcal N^{(1)}_ M( P)$ is defined as the set of all $x\in  M$ for which there exist $x_1,x_2,...,x_k\in  M$ such that  $ P x\subseteq \sum_i x_i  P$. 

We begin this subsection with two general results on quasi-normalizers. 
Firstly, we record the following formula for one sided quasi-normalizers of corners.

\begin{lem}[\!\!\cite{Po03, FGS10}]\label{quasi}
Let $ P\subseteq M$ be an inclusion of tracial von Neumann algebras. 

Then $\text{W}^*(q\mathcal N^{(1)}_{p M p}(p P p))=p\;\text{W}^*(q\mathcal N^{(1)}_{ M}( P))\;p$, for any projection $p\in P$.

Moreover, $q\mathcal N_{p' M p'}^{(1)}( P p')=p'\;q\mathcal N_{ M}^{(1)}( P)\;p'$, for any projection $p'\in P'\cap  M$.
\end{lem}

\noindent
The main assertion follows from the proof of \cite[Lemma 3.5]{Po03}, where a similar formula for the usual quasi-normalizer is provided.
It appears as such in \cite[Proposition 6.2]{FGS10}. The moreover assertion is immediate.

Secondly, we establish a useful property of subalgebras having a trivial one sided quasi-normalizer.

\begin{lem}\label{cornertowhole}
Let $ M$ be a tracial von Neumann algebra, and $ P, Q\subseteq  M$ von Neumann subalgebras.
Assume that  $q\mathcal N_{ M}^{(1)}( P)= P$ and $ Q$ is a II$_1$ factor. Suppose also that  $ P\prec_{ M}^{s} Q$
and that $q P q=q Q q$, for some non-zero projection $q\in P$.

 Then there exists $u\in\mathcal U( M)$ such that $u P u^*= Q $.
 Moreover, if $ P\subseteq Q$, then $ P= Q$.
\end{lem}

{\it Proof.}
Let us first show that $ P$ can be unitarily conjugated into $ Q$. To this end, let $r\in\mathcal Z( P)$ be a non-zero projection. Since $ P r\prec_{ M} Q$, we can find projections $r_0\in  P r$ and $q_0\in Q$, a non-zero partial isometry $v\in q_0 M r_0$, and a $*$-homomorphism $\theta:r_0 P r_0\rightarrow q_0 Q q_0$ such that $\theta(x)v=vx$, for all $x\in r_0 P r_0$. 
Moreover, after replacing $r_0$ with a smaller projection, we may assume that $\uptau(q_0)\leqslant\uptau(q)$. Since $ Q$ is a II$_1$ factor we can find a unitary $\eta\in Q$ such that $q_1:=\eta q_0\eta^*\leqslant q$. Let $\varphi:r_0 P r_0\rightarrow q_1 Q q_1$ be given by $\varphi(x)=\eta \theta(x)\eta^*$. If we put 
$w=\eta v\in q_1 M r_0$, then $\varphi(x)w=wx$, for all $x\in r_0 P r_0$.
Since $q_1\leqslant q$ we have that $q_1 Q q_1=q_1 P q_1$, and thus $w P r_0=\varphi(p_0Pr_0)w\subseteq  P w$. 

We claim that $w\in P$.
This follows from the proof of \cite[Lemma 3.5]{Po03}. For completeness, we include the argument here.
Let $z\in  \mathcal Z(P)$ be the central support of $r_0$.  If $\varepsilon>0$, then we can find a projection $z'\in\mathcal Z( P)z$ such that $\uptau(z-z')<\varepsilon$ and there exists partial isometries $\xi_1,...,\xi_n\in r_0 P$ satisfying $\sum_{i\geq 1}\xi_i^*\xi_i=z'$. Then $wz' P=w P z'\subseteq\sum_{i=1}^nw P\xi_i^*\xi_i\subseteq\sum_{i=1} P w\xi_i,$ and therefore $wz'\in q\mathcal N_{ M}^{(1)}( P)= P$. Since $\varepsilon>0$ is arbitrary and $w=wz$, the claim follows.

Put $r_1=w^*w\in r_0 P r_0$ and $q_2=ww^*\in q_1 P q_1=q_1 Q q_1$. Then $w P w^*=q_2 P q_2=q_2 Q q_2$, so in particular $r_1 P r_1$ can be unitarily conjugated into $ Q$. Since $ Q$ is a II$_1$ factor,  we deduce that $ P r'$ can be unitarily conjugated into $ Q$, where $r'\in \mathcal Z( P)r$ denotes the central support of $r_1$. Thus, for every non-zero projection $r\in\mathcal Z( P)$, there is a non-zero projection $r'\in\mathcal Z( P)r$ such that $ P r'$ can be unitarily conjugated into $ Q$. Since $ Q$ is a II$_1$ factor, a maximality argument implies the existence of $u\in\mathcal U(M)$ such that $u P u^*\subseteq Q$.

Finally, we prove that $u P u^*= Q$, which will imply both assertions of the lemma. Let $r_0\in P$ be a projection such that $\uptau(r_0)\leqslant\uptau(q)$ and put $q_0=ur_0u^*\in Q$. Then we can find $\eta\in\mathcal U( Q)$  such that $q_1:=\eta q_0\eta^*\leqslant q$. Since $ur_0 P r_0u^*\subseteq q_0 Q q_0= \eta^* q_1 Q q_1\eta=\eta^*q_1 P q_1\eta$, it follows as above that $\eta ur_0\in P$. This implies that $ur_0 P r_0u^*=\eta^*q_1 P q_1\eta$, thus $ur_0 P r_0u^*=q_0 Q q_0$. Hence we have that $r_0 P r_0=r_0(u^* Q u)r_0$, for any projection $r_0\in P$ with $\uptau(r_0)\leqslant\uptau(q)$. This clearly implies that $ P=u^* Q u$, which finishes the proof.
\hfill$\blacksquare$

In the rest of this subsection, we establish several results controlling quasi-normalizers in group von Neumann algebras.

\begin{lem}[\!\!\cite{Po03}]\label{mixing}
Let $\Gamma_1<\Gamma$ be countable groups, and $ P\subseteq pL(\Gamma_1)p$ be a von Neumann subalgebra, for a projection $p\in L(\Gamma_1)$.
Assume that $ P\nprec_{L(\Gamma_1)}L(\Gamma_1\cap g\Gamma_1 g^{-1})$, for all $g\in\Gamma\setminus\Gamma_1$.

If
 $x\in L(\Gamma)$ satisfies $x P\subseteq \sum_{i=1}^nL(\Gamma_1)x_i$,  for some $x_1,...,x_n\in L(\Gamma)$, then $xp\in L(\Gamma_1)$. 
\end{lem}

{\it Proof.}
Since  $ P\nprec_{L(\Gamma_1)}L(\Gamma_1\cap g\Gamma_1 g^{-1})$, for all $g\in\Gamma\setminus\Gamma_1$, we can find a net $u_n\in\mathcal U(P)$ such that 
$\|E_{L(\Gamma_1\cap g\Gamma_1 g^{-1})}(u_na)\|_2\rightarrow 0$, for all $a\in L(\Gamma_1)$ and $g\in\Gamma\setminus\Gamma_1$ 
(see the proof \cite[Theorem 4.3]{IPP05} and also \cite[Remark 2.3]{DHI16}).
By a result of Popa (see \cite[Theorem 3.1]{Po03} and also \cite[Theorem 1.1]{IPP05} and \cite[Lemma D.3]{Va06}), in order to get the conclusion it suffices to show that
\begin{equation}\label{weakly}\|E_{L(\Gamma_1)}(bu_nc)\|_2\rightarrow 0,\;\;\;\text{for every $b,c\in L(\Gamma)$ with $E_{L(\Gamma_1)}(b)=0$}.
\end{equation}
\noindent
By Kaplansky's density theorem, in order to prove \eqref{weakly}, we may assume that $b=u_g, c=u_h$, for some $g\in\Gamma\setminus\Gamma_1$ and $h\in\Gamma$. 
If $g\Gamma_1 h\cap\Gamma_1=\emptyset$, then $E_{L(\Gamma_1)}(u_gu_nu_h)=0$, for all $n$. If $g\Gamma_1 h\cap\Gamma_1$ is non-empty, fix $k\in g\Gamma_1 h\cap\Gamma_1$, and put $l=g^{-1}kh^{-1}\in\Gamma_1$. 
Then $g\Gamma_1h\cap\Gamma_1=(g\Gamma_1 g^{-1}\cap\Gamma_1)k$. Thus, if $\gamma\in\Gamma_1$, then $g\gamma h\in g\Gamma_1 h\cap\Gamma_1$ if and only if 
$\gamma\in (\Gamma_1\cap g^{-1}\Gamma_1g)l$. Therefore, $$\|E_{L(\Gamma_1)}(u_gu_nu_h)\|_2=\|E_{L(\Gamma_1\cap g^{-1}\Gamma_1g)}(u_nu_l^*)\|_2\rightarrow 0.$$
\noindent
This altogether proves \eqref{weakly}, and finishes the proof.
\hfill$\blacksquare$

\noindent
The next result strengthens the conclusion of Lemma \ref{mixing} in the case of inclusions of group von Neumann algebras. 

\begin{lem}\label{mixing2}
Let $\Gamma_1<\Gamma_2<\Gamma$ be countable groups. Denote by $S\subseteq\Gamma$ the set of $g\in\Gamma$ such that $[\Gamma_1:\Gamma_1\cap g\Gamma_2g^{-1}]<\infty$.

If  $x\in L(\Gamma)$ satisfies $L(\Gamma_1)x\subseteq\sum_{i=1}^nx_iL(\Gamma_2)$, for some $x_1,...,x_n\in L(\Gamma)$, then $x$ belongs to the $\|.\|_2$-closure of the linear span of $\{u_g\}_{g\in S}$.
\end{lem}

This result generalizes \cite[Theorem 5.1]{FGS10}, which addressed the case $\Gamma_1=\Gamma_2$. Although later on we will only use this particular case of Lemma \ref{mixing2}, for completeness we provide a different proof that at the same time handles the general case.  The proof that we include follows closely the proof of \cite[Theorem 2.1]{Po03}.

{\it Proof.}
Let $\mathcal K\subseteq\ell^2(\Gamma)$ be the $\|.\|_2$-closure of the linear span of $L(\Gamma_1)xL(\Gamma_2)$, and $f$ the orthogonal projection from $\ell^2(\Gamma)$ onto $\mathcal K$. Since $\mathcal K$ is a $L(\Gamma_1)$-$L(\Gamma_2)$-bimodule, $f\in L(\Gamma_1)'\cap \langle L(\Gamma),e_{L(\Gamma_2)}\rangle$.  Since $\mathcal K$ is contained in the $\|.\|_2$-closure of $\sum_{i=1}^nx_iL(\Gamma_2)$, we also have that $Tr(f)<\infty$. Viewing $f$ as an element of $L^2(\langle L(\Gamma_1),e_{L(\Gamma_2)}\rangle)$, Claim \eqref{invariant} gives that $f$ belongs to the $\|.\|_{2,Tr}$-closure of the linear span of $\{u_geu_h\}_{g\in S, h\in\Gamma}$. This implies that $f(\ell^2(\Gamma))$ is contained  the $\|.\|_2$-closure of the linear span of $\{u_g\}_{g\in S}$. Since $x\in f(\ell^2(\Gamma))$, the conclusion follows.
\hfill$\blacksquare$

\begin{cor}[\!\!\cite{FGS10}] \label{FGS}Let $\Sigma<\Gamma$ be countable groups.

Then $W^*(q\mathcal N^{(1)}_{L(\Gamma)}(L(\Sigma)))=L(\Delta)$, where $\Delta<\Gamma$ is the subgroup generated by $\text{Comm}_{\G}^{(1)}(\Sigma)$.

In particular, if $\text{Comm}_{\G}^{(1)}(\Sigma)=\Sigma$, then $q\mathcal N^{(1)}_{L(\Gamma)}(L(\Sigma))=L(\Sigma)$.
\end{cor}

{\it Proof.}
This result is part (ii) of \cite[Corollary 5.2]{FGS10}. 
For completeness, we show how  it follows from Lemma \ref{mixing2}. 
If $g\in \text{Comm}^{(1)}_{\Gamma}(\Sigma)$, then $u_g\in q\mathcal N^{(1)}_{L(\Gamma)}(L(\Sigma))$. This implies the inclusion $\supseteq$. If $x\in q\mathcal N^{(1)}_{L(\Gamma)}(L(\Sigma))$, then Lemma \ref{mixing2} gives that $x\in L(\Delta)$, which implies the reverse inclusion. 
\hfill$\blacksquare$

We end this section with two results concerning von Neumann algebras of amalgamated free product groups.

\begin{cor}[\!\!\cite{IPP05}]\label{mixing3}
Let $\Gamma=\Gamma_1\ast_{\Sigma}\Gamma_2$ be an amalgamated free product group. Let $ P\subseteq p(\Gamma_1)p$ be a von Neumann subalgebra, for a projection $p\in P$.
 Assume that $ P\nprec_{L(\Gamma_1)}L(\Sigma)$.

 If $x\in L(\Gamma)p$ satisfies $x P\subseteq\sum_{i=1}^nL(\Gamma_1)x_i$, for some $x_1,...,x_n\in L(\G)$, then $x\in L(\Gamma_1)$.

\end{cor}

{\it Proof.}
This result is a particular case of \cite[Theorem 1.1]{IPP05}. Since $\Gamma_1\cap g\Gamma_1 g^{-1}\subseteq\Sigma$, for all $g\in\Gamma\setminus\Gamma_1$, it also follows from Lemma \ref{mixing}.  
\hfill$\blacksquare$

\begin{lem}\label{mixing4}
Let $\Gamma=\Gamma_1\ast_{\Sigma}\Gamma_2$ be an amalgamated free product group with $\text{Comm}^{(1)}_{\G_i}(\Sigma)=\Sigma$, for  every $i\in\overline{1,2}$. 
Then we have the following:
\begin{enumerate}
\item $\text{Comm}^{(1)}_{\G}(\Sigma)=\Sigma$.
\item $L(\Sigma)\nprec_{L(\G)}L(\Sigma\cap g\Sigma g^{-1})$, for every $g\in\Gamma\setminus\Sigma$.
\item $L(\Sigma)\nprec_{L(\G)}L(\Gamma_i\cap g\Gamma_ig^{-1})$, for every $g\in\Gamma_i\setminus\Sigma$ and $i\in\overline{1,2}$.
\end{enumerate}
\end{lem}

{\it Proof.}
(1) Let $g\in\Gamma\setminus\Sigma$. Let $g=g_1g_2...g_n$ be the reduced form of $g$, where $n\geq 1$ is an integer, $j(k)\in\overline{1,2}$ and $g_k\in \Gamma_{j(k)}\setminus\Sigma$,  for all $k\in\overline{1,n}$, and $j(1)\not=j(2)\not=...\not=j(n)$.
If $x\in\Sigma\cap g\Sigma g^{-1}$, then $g^{-1}xg=g_n^{-1}...g_2^{-1}g_1^{-1}xg_1g_2...g_n\in\Sigma$, which forces $g_1^{-1}xg_1\in\Sigma$. Thus, $\Sigma\cap g\Sigma g^{-1}\subseteq \Sigma\cap g_1\Sigma g_1^{-1}$. Since $g\in\Gamma_{j(1)}\setminus\Sigma$, we have that $[\Sigma:\Sigma\cap g_1\Sigma g_1^{-1}]=\infty$, which implies that $g\not\in \text{Comm}^{(1)}_{\G}(\Sigma)$.

(2) Let $g\in\Gamma$ such that $L(\Sigma)\prec_{L(\G)}L(\Sigma\cap g\Sigma g^{-1})$. By applying Lemma \ref{intertwiningroups}, we find $h\in\Gamma$ such that $[\Sigma:\Sigma\cap h(\Sigma\cap g\Sigma g^{-1})h^{-1}]<\infty$, and thus $[\Sigma:\Sigma\cap h\Sigma h^{-1}]<\infty$ and $[\Sigma:\Sigma\cap (hg)\Sigma (hg)^{-1}]<\infty$. 
By  using part (1) we deduce that $h,hg\in\Sigma$, and thus  $g\in\Sigma$.

(3) Assume that $L(\Sigma)\prec_{L(\G)}L(\Gamma_i\cap g\Gamma_ig^{-1})$, for some $g\in\Gamma\setminus\Sigma$ and $i\in\overline{1,2}$. 
Let $g=g_1g_2...g_n$ be the reduced form of $g$, where $n\geq 1$ is an integer, $j(k)\in\overline{1,2}$ and $g_k\in \Gamma_{j(k)}\setminus\Sigma$,  for all $k\in\overline{1,n}$, and $j(1)\not=j(2)\not=...\not=j(n)$.
  Let $x\in \Gamma_i\cap g\Gamma_i g^{-1}$. Then $g^{-1}xg=g_n^{-1}...g_2^{-1}g_1^{-1}xg_1g_2...g_n\in\Gamma_i$. Since $g\in\Gamma\setminus\Gamma_i$, we can find $k\in\overline{1,n}$ with $j(k)\not=i$. If $j(1)=i$, then we must have that $g_1^{-1}xg_1\in\Sigma$ and $g_2^{-1}g_1^{-1}xg_1g_2\in\Sigma$, and hence $\Gamma_i\cap g\Gamma_ig^{-1}\subseteq g_1(\Sigma\cap g_2\Sigma g_2^{-1})g_1^{-1}$. If $j(1)\not=i$, then we must have that $x\in\Sigma$ and $g_1^{-1}xg_1\in\Sigma$, and hence $\Gamma_i\cap g\Gamma_ig^{-1}\subseteq \Sigma\cap g_1\Sigma g_1^{-1}$. In either case, we would conclude that $L(\Sigma)\prec_{ M}L(\Sigma\cap h\Sigma h^{-1})$, for some $h\in\Gamma\setminus\Sigma$. By (2) this is a contradiction.
\hfill$\blacksquare$

\subsection{Almost malnormality of maximal amenable subgroups in hyperbolic groups}\label{exb}
In this section, we justify an assertion made in Example \ref{ex}(b). 
\begin{lem}\label{exemplu}
Let  $C<D$ be an infinite maximal amenable subgroup of a hyperbolic group $D$.

Then $C\cap gCg^{-1}$ is finite, for every $g\in D\setminus C$. 
\end{lem}
This result is likely well-known, but for lack of a reference, we include a proof.
Note that it implies that if $\Sigma_0:=A\wr C$ and $\Gamma_0:=A\wr D$,
then $[\Sigma_0:\Sigma_0\cap g\Sigma_0 g^{-1}]=\infty$, for all $g\in\Gamma_0\setminus\Sigma_0$, as claimed in Example \ref{ex}(b).

{\it Proof.}
By \cite[Th\'{e}or\`{e}me 8.37]{GdH90}, $C$ admits an infinite cyclic subgroup $C_0=\{a^n|n\in\mathbb Z\}$  of finite index. By \cite[Th\'{e}or\`{e}me 8.29]{GdH90}, if $\partial D$ denotes the boundary of $D$, then $\overline{C_0}\cap\partial D$ contains exactly two points $\{x_1,x_2\}$, both fixed by $a$. 

We claim that if $g\in D$ and $C_0\cap gC_0g^{-1}\not=\{e\}$, then $g$ stabilizes the set $\{x_1,x_2\}$.
Let $m,n\in\mathbb Z\setminus \{0\}$ such that $ga^mg^{-1}=a^n$.
If $i\in\{1,2\}$, we can find a sequence $\{p_k\}$ such that $x_i=\lim\limits_{k\ra\infty}a^{p_k}$. Then $x_i=\lim\limits_{k\ra\infty}a^{m\lfloor\frac{p_k}{m}\rfloor}$ and thus $gx_i=\lim\limits_{k\ra\infty}ga^{m\lfloor\frac{p_k}{m}\rfloor}=\lim\limits_{k\ra\infty}a^{n\lfloor\frac{p_k}{m}\rfloor}g=\lim\limits_{k\ra\infty}a^{n\lfloor\frac{p_k}{m}\rfloor}\in\{x_1,x_2\}$. 

Since $[C:C_0]<\infty$, the claim implies that $C\subseteq\text{Stab}_D\{x_1,x_2\}$.  
By \cite[Th\'{e}or\`{e}me 8.30]{GdH90}, $\text{Stab}_D\{x_1,x_2\}$ contains a finite index cyclic subgroup, hence is amenable.
Thus, we get that $C=\text{Stab}_D\{x_1,x_2\}$. Using the claim again gives that  $C\cap gCg^{-1}$ is finite, for all $g\in D\setminus C$.
\hfill$\blacksquare$

\subsection{Property (T) and Haagerup's property for groups and algebras} In this section, we record the well-known relationship between property (T) and Haagerup's property for countable groups and their 
von Neumann algebras. We refer the reader to \cite{Po01} for the definitions of these notions. As shown in \cite{CJ85} and \cite{Ch83} a countable icc group $\Gamma$ has property (T) and respectively Haagerup's property if and only if the II$_1$ factor $L(\Gamma)$ does. Moreover, this result holds if $\Gamma$ is not necessarily icc (see \cite[Propositions 3.1 and 5.1]{Po01}).
Here we note that arguments from \cite{Po01} show that the result remains true if in addition $L(\Gamma)$ is replaced by one of its corners.

\begin{lem}\label{HT} Let $\Gamma$ be a countable group and $p\in L(\Gamma)$ be a non-zero projection. Then \begin{enumerate} \item $\Gamma$ has property (T) if and only if $pL(\Gamma)p$ does. \item $\Gamma$ has Haagerup's property if and only if $pL(\Gamma)p$ does. \end{enumerate} \end{lem}

{\it Proof.}
(1) Assume that $\Gamma$ has property (T). Then \cite[Proposition 5.1]{Po01} implies that $L(\Gamma)$ has property (T), and \cite[Proposition 4.7 (2)]{Po01} further gives that $pL(\Gamma)p$ has property (T).

\noindent
Conversely, assume that $pL(\Gamma)p$ has property (T). Denoting by $z\in L(\Gamma)$ the central support of $p$, \cite[Proposition 4.7 (3)]{Po01} implies that $L(\Gamma)z$ has property (T).
Let $\varphi_n:\Gamma\rightarrow\mathbb C$ be a sequence of positive definite functions  such that $\varphi_n(e)=1$, for all $n$, and $\varphi_n(g)\rightarrow 1$, for all $g\in\Gamma$. Then the formula $\Phi_n(x)=\sum_g\varphi_n(g)a_gu_g$, for every $x=\sum_ga_gu_g\in L(\G)$, defines a sequence $\Phi_n:L(\Gamma)\rightarrow L(\Gamma)$ of unital, tracial, completely positive maps such that $\|\Phi_n(x)-x\|_2\rightarrow 0$, for every $x\in L(\Gamma)$. 
Thus, if we let $\Psi_n(x)=\Phi_n(x)z$, for every $x\in L(\G)z$, then $\Psi_n:L(\Gamma)z\rightarrow L(\Gamma)z$ is a sequence of unital, subtracial, completely positive maps such that $\|\Psi_n(x)-x\|_2\rightarrow 0$, for every $x\in L(\G)z$. Since $L(\G)z$ has property $(T)$, we get that $\sup \{\|\Psi_n(u)-u\|_2\;|\;u\in\mathcal U(L(\G)z)\}\rightarrow 0$. Since $\sup \{\|\Phi_n(uz)-\Phi_n(u)z\|_2\;|\;u\in\mathcal U(L(\G))\}\rightarrow 0$, we get that $\sup \{\|\Phi_n(u_g)z-u_gz\|_2\;|\;g\in\Gamma\}\rightarrow 0$. As $\|\Phi_n(u_g)z-u_gz\|_2=|\varphi_n(g)-1|\;\|z\|_2$, for all $g\in\Gamma$, we conclude that $\sup\{|\varphi_n(g)-1|\;|\;g\in\Gamma\}\rightarrow 0$. Thus, $\Gamma$ has property (T).

\noindent
(2) Assume that $\Gamma$ has Haagerup's property.  Then \cite[Propositions 3.1 and 2.4 (1)]{Po01} together imply that $pL(\G)p$ has Haagerup's property.

\noindent
Conversely, assume that $pL(\G)p$ has Haagerup's property.  Denoting by $z\in L(\Gamma)$ the central support of $p$, \cite[Proposition 2.4 (2)]{Po01} implies that $L(\Gamma)z$ has Haagerup's property. Let $\Phi_n:L(\Gamma)z\rightarrow L(\Gamma)z$ be a sequence of completely positive maps such that $\tau\circ\Phi_n\leqslant\tau$ and the set $\{\Phi_n(x)\;|\;x\in L(\Gamma)z, \|x\|\leqslant 1\}$ is $\|.\|_2$-precompact, for all $n$, and $\|\Phi_n(x)-x\|_2\rightarrow 0$, for all $x\in L(\G)z$. Then $\varphi_n:\Gamma\rightarrow\mathbb C$ given by $\varphi_n(g)=\tau(z)^{-1}\tau(\Phi_n(u_gz)(u_gz)^*)$ is a $c_0$ positive definite function. As $\varphi_n(g)\rightarrow 1$, for all $g\in\Gamma$, we get that  $\Gamma$ has Haagerup's property.  
\hfill$\blacksquare$


\section{Identification of Peripheral Subgroups via $W^*$-Equivalence}
\noindent
In this section we establish the main technical result needed in the proof of Theorem \ref{A}.
Throughout the section, we will work with amalgamated free product groups $\Gamma$ satisfying the following:

\begin{assumption}\label{AFP}
 $\Gamma=\Gamma_1\ast_{\Sigma}\Gamma_2$ is an amalgamated free product group, where
\begin{enumerate} 
\item $\Sigma$ is a common amenable subgroup of $\Gamma_1$ and $\Gamma_2$.
\item $\Gamma_i=\Gamma_i^1\times\Gamma_i^2$, where $\Gamma_i^j$ is an icc, non-amenable, bi-exact group, for every $i,j\in\overline{1,2}$. 
\end{enumerate}
\end{assumption}

\noindent
The main goal of this section to establish the following structural result for groups $\Lambda$ in the $W^*$-equivalence class of an amalgamated free product group $\Gamma$ as in Assumption \ref{AFP}. 

 
\begin{theorem}\label{peripheralidentification2}
Let  $\G=\G_1\ast_{\Sigma}\G_2$ be as in Assumption \ref{AFP}, and put $ M=L(\G)$. Let $\La$ be an arbitrary group such that $ M=L(\Lambda)$. 
Then one of the following two conditions holds:

\begin{enumerate}
\item For every $i\in\overline{1,2}$, there exists a subgroup $\Th_i<\La$ such that $u_i L(\Th_i)u^*_i=L(\Gamma_i)$ for some $u_i\in\mathcal U( M)$.
\item There exist a subgroup $\Th<\La$, non-zero projections $r_1,r_2\in\mathcal Z(L(\Th))$ with $r_1+r_2=1$, and $u\in\mathcal U( M)$ such that $q_i:=ur_iu^*\in L(\Gamma_i)$ and $uL(\Th)r_iu^*=q_iL(\G_i)q_i$, for every $i\in\overline{1,2}$.
\end{enumerate} 
\end{theorem}

\noindent
The rest of this section is devoted to the proof of Theorem \ref{peripheralidentification2}.
 The starting point is the following key result showing that any such group $\La$ must contain commuting non-amenable subgroups.

\begin{theorem}\label{comm1} Let  $\G=\G_1\ast_{\Sigma}\G_2$ be as in Assumption \ref{AFP}. Put $ M=L(\G)$, and let $\La$ be an arbitrary group such that $ M=L(\Lambda)$. 

Then for every $i, j\in \overline{1,2}$  we can find a non-amenable subgroup $\Delta <\La$  such that  $C_{\La}(\Delta)$ is non-amenable and $L(\Gamma_i^{j}) \prec_ M L(\Delta)$.  
 \end{theorem}

\noindent
There are two main ingredients in the proof of Theorem \ref{comm1}. 
Thus, we first use repeatedly Ozawa's  work on the structure of subalgebras with non-amenable commutant inside von Neumann algebras of  relatively bi-exact groups \cite{Oz03,Oz04, BO08} (see also the more recent developments \cite{CS11, CSU11}).  
We refer the reader to \cite[Definition 15.1.2]{BO08} for the notion of relative bi-exactness  for groups.  

 A second crucial ingredient is the ultrapower technique for group von Neumann algebras introduced by the second author in \cite[Theorem 3.1]{Io11}. Note that this technique has recently been used in several other works \cite{CdSS15,KV16,DHI16}.

\subsection*{Proof of Theorem \ref{comm1}.}
Denote by $\set{u_g}_{g\in \Gamma}$ and $\set{v_h}_{h\in \Lambda} $ the canonical unitaries generating $ M$.  Following \cite{PV09}, we consider a $*$-homomorphism $ \triangle :  M\to  M\bar{\otimes}  M$, called the comultiplication along $\Lambda$, and defined by $\Delta(v_h)=v_h\otimes v_h$, for all $h\in\La$. Then we have

\begin{claim}\label{2001} For every $i,j\in\overline{1,2}$, there exist $m,n\in\overline{1,2}$  such that $\triangle(L(\Gamma_i^j))\prec_{ M\bar\otimes  M}  M\bar\otimes L(\Gamma_m^n)$.
\end{claim}

\noindent \emph{Proof of Claim \ref{2001}.} Let $i\in\overline{1,2}$, and denote $ P=L(\Gamma_i^1)$, $ Q=L(\Gamma_i^2)$. 
Then $\triangle( P)$ and $\triangle( Q)$ are commuting non-amenable subalgebras of ${ M\bar\otimes  M}=L(\G\times \G)$. On the other hand,  by \cite[Lemma 15.3.3 and Proposition 15.3.12]{BO08}, $\Gamma\times\Gamma$ is bi-exact relative to the family $\mathcal G = \{\G\times \Gamma_m|m\in \overline{1,2}\}\cup\{ \Gamma_m\times \G|m \in \overline{1,2}\}.$  

 By applying \cite[Theorem 15.1.5]{BO08}, we deduce that $\triangle( P) \prec_{ M\bar\otimes  M} L(G)$, for some  $G\in\mathcal G$.
 Since the flip automorphism of $ M\bar{\otimes} M$ acts identically on $\Delta( P)$, we may assume that $G=\G\times \Gamma_m$, for $m\in\overline{1,2}$.  
 Thus, there exist projections $p\in \triangle( P), q\in  M\bar{\otimes}L(\Gamma_m)$, a non-zero partial isometry $v\in q( M\bar\otimes  M)p$, and a $\ast$-homomorphism $\varphi:p\triangle( P)p\ra q( M\bar{\otimes}L(\Gamma_m))q$ such that 
 \begin{equation}\label{1002}\varphi(x)v=vx,\text{ for all }x\in p\triangle( P)p.\end{equation} Notice that $vv^*\in\varphi(p\triangle ( P)p)'\cap q( M\bar\otimes  M)q$ and $v^*v\in (p\triangle( P)p)'\cap p( M\bar \otimes  M)p$. 
 We may assume that $q$ is equal to the support projection of $E_{q( M\bar{\otimes}L(\Gamma_m))q}(vv^*)$.
Let us show that
 \begin{equation}\label{1001}\varphi(p\triangle ( P)p)\nprec_{ M\bar{\otimes}L(\Gamma_m)}  M\bar{\otimes} L(\Sigma).\end{equation} Indeed, if (\ref{1001}) does not hold, \cite[Lemma 1.12]{IPP05}  would imply that $\triangle( P) \prec_{ M\bar{\otimes} M}  M\bar{\otimes}L(\Sigma)$.  Since $ P$ has no amenable direct summand, by \cite[Proposition 7.2(4)]{IPV10}, this would contradict the amenability of $\Sigma$.

Since $ M\bar{\otimes} M=L(\Gamma\times\Gamma_1)*_{L(\Gamma\times\Sigma)}L(\Gamma\times\Gamma_2)$, by combining \eqref{1001} and Corollary \ref{mixing3} we conclude that $\varphi(p\triangle ( P)p)'\cap q( M\bar{\otimes} M)q\subset q( M\bar{\otimes}L(\Gamma_m))q$, hence $vv^*\in  M\bar{\otimes}L(\Gamma_m)$.
Thus, equation (\ref{1002}) implies that $v\triangle( P)v^*\subseteq  M\bar{\otimes}L(\Gamma_m)$. 
Since $\triangle( P) \nprec_{ M\bar{\otimes} M}  M\bar{\otimes}L(\Sigma)$, applying Corollary \ref{mixing3} again gives that $v(\triangle( P)\vee (\triangle( P)'\cap  M\bar{\otimes} M))v^*\subset M\bar{\otimes}L(\Gamma_m)$. 
Since $ M\bar{\otimes}L(\Gamma_m)$ is a factor, we can thus find a non-zero projection $e\in \mathcal Z(\triangle( P)'\cap ( M\bar \otimes  M) )$  and $u\in\mathcal U( M\bar{\otimes} M)$ such that \[u(\triangle( P)\vee (\triangle ( P)'\cap  M\bar \otimes  M) )eu^*\subseteq  M\bar{\otimes}L(\Gamma_m).\]

In particular, we get that 
$u\triangle( P\vee  Q)eu^*\subseteq  M\bar{\otimes}L(\Gamma_m).$
Thus, $u\triangle( P)eu^*$ and $u\triangle( Q)eu^*$ are commuting non-amenable subfactors of $ M\bar{\otimes}L(\Gamma_m)=L(\Gamma\times\Gamma_m)$. 
Since $\Gamma\times\Gamma_m$ is bi-exact relative to the family $\mathcal H=\{\G\times \Gamma_m^1,\Gamma\times \Gamma_m^2, \Gamma_1\times \Gamma_m,\Gamma_2\times\Gamma_m\}$, by applying  \cite[Theorem 15.1.5]{BO08}  we get that $\triangle( P)e\prec_{ M\bar \otimes  M} L(H)$,  for some $H\in \mathcal H.$

If $H=\G\times \Gamma_m^1$ or $H=\Gamma\times \Gamma_m^2$, then the claim follows. Therefore, it remains to analyze the case when $H=\Gamma_r\times \Gamma_m$, for some $r\in \overline{1,2}$. In this case, by arguing as above, we find  a non-zero projection $f\in \mathcal Z((\Delta( P)e)'\cap e( M\bar \otimes  M) e) $  and $w\in \mathcal U( M\bar \otimes  M)$ such that we have
$w\triangle( P\vee  Q) f w^*\subseteq L(\Gamma_r)\bar{\otimes}L(\Gamma_m).$  In particular, $w\triangle( P)fw^*$ and $w\triangle( Q)fw^*$ are commuting, non-amenable subfactors  of $L(\Gamma_r\times\Gamma_m)$

Since $\Gamma_r\times\Gamma_m$ is bi-exact relative to  $\mathcal K=\{\Gamma_r^1\times\Gamma_m,\Gamma_r^2\times\Gamma_m,\Gamma_r\times\Gamma_m^1,\Gamma_r\times\Gamma_m^2 \}$, \cite[Theorem 15.1.5]{BO08} implies that $\triangle( P)f\prec_{ M\bar \otimes  M} L(K)$ and hence $\triangle( P)\prec_{ M\bar\otimes M} L(K)$, for some $K\in \mathcal K.$
Since the flip automorphism of $ M\bar{\otimes} M$ acts identically on $\Delta( M)$, this concludes the proof of the claim.$\hfill\square$

We are now in position to apply the ultrapower technique from \cite{Io11}, which we recall in the following form. 
This result is essentially contained in the proof of \cite[Theorem 3.1]{Io11}. Stated as such, it is a particular case of \cite[Theorem 4.1]{DHI16}.

\begin{theorem}[\!\!\cite{Io11}]\label{Io11}
Let $\Lambda$ be a countable group, $ M=L(\Lambda)$, and $\{v_h\}_{h\in\La}$ the canonical unitaries generating $M$.
Let  $\triangle: M\rightarrow M\bar{\otimes} M$ be the $*$-homomorphism given by $\triangle(v_h)=v_h\otimes v_h$, for all $h\in\La$. Let $A, B\subseteq M$ be von Neumann subalgebras such that $\triangle(A)\prec M\bar{\otimes}B$.

Then there exists a decreasing sequence of subgroups $\La_k<\Lambda$ such that $A\prec_{ M}L(\La_k)$, for every $k\geq 1$, and $B'\cap M\prec_{ M}L(\cup_{k\geq 1}C_{\Lambda}(\La_k))$.
\end{theorem}

 Going back to the proof of Theorem \ref{comm1}, by combining Claim (\ref{2001}) and Theorem \ref{Io11}, we deduce the existence of a decreasing sequence of subgroups $\La_k<\La$ such that $L(\Gamma_i^j)\prec_{ M}L(\La_k)$, for every $k\geq 1$, and $L(\Gamma_m^n)'\cap M\prec_{ M}L(\cup_{k\geq 1}C_{\Lambda}(\La_k))$.
Since $L(\Gamma_m^n)'\cap M$ is non-amenable, as it contains $L(\Gamma_m^s)$, where $\{n,s\}=\{1,2\}$, we get that $\cup_{k\geq 1}C_{\Lambda}(\La_k)$ is non-amenable. Thus, there is $k\geq 1$ such that $C_{\Lambda}(\La_k)$ is non-amenable. Letting $\Delta:=\La_k$ the conclusion follows since $L(\Gamma_i^j)\prec_{ M}L(\Delta)$ and in particular $\Delta$ is non-amenable.
 \hfill$\blacksquare$
  
 We continue with the second step towards proving Theorem \ref{peripheralidentification2}. More precisely, we use the commuting subgroups of the mysterious group $\La$ provided by Theorem \ref{comm1}, to identify the algebras of the peripheral subgroups $\Gamma_1,\Gamma_2$ of  $\G$ with algebras of certain subgroups of $\La$. Our proof is inspired by the analysis performed in \cite[Theorem 4.3]{CdSS15}.

\begin{theorem}\label{peripheralidentification1}
Let  $\G=\G_1\ast_{\Sigma}\G_2$ be as in Assumption \ref{AFP}, and put $ M=L(\G)$. Let $\La$ be a group such that $ M=L(\Lambda)$.  Assume that there exists a non-amenable subgroups $\Delta<\La$ such that $C_{\La}(\Delta)$ is non-amenable and $L(\Gamma^1_1)\prec_{ M} L(\Delta)$. 

Then we can find a 
group $\Th<\La$, 
 projections $ r_1,r_2\in \mathcal Z (L(\Theta))$ with $r_1\not=0$, $r_1+r_2=1$, and $u\in\mathcal U( M)$ such that  $ur_1u^*\in L(\Gamma_1)$, $uL(\Th)r_1u^*=ur_1u^*L(\Gamma_1)ur_1u^*$ and $uL(\Th)r_2u^*\subseteq L(\Gamma_2)$.

\end{theorem}

{\it Proof.} Let $\Omega$ be the group of $h\in\La$ such that $\{\delta h\delta^{-1}|\delta\in\Delta\}$ is finite. Then $\Omega$ is normalized by $\Delta$, hence $\Delta\Omega$ is a subgroup of $\La$.
Let $\Theta<\Lambda$ be the subgroup generated by $\text{Comm}^{(1)}_{\La}(\Delta\Omega)$.
By Corollary \ref{FGS} we have that \begin{equation}\label{theta} \text{W}^*(q\mathcal N^{(1)}_ M(L(\Delta\Omega)))=L(\Theta).\end{equation}
\noindent
Since $L(\Gamma_2)$ is a II$_1$ factor, there is a maximal projection $r_2\in\mathcal Z (L(\Theta))$ such that  $L(\Theta)r_2$ can be unitarily conjugated into $L(\Gamma_2)$.  Let $w_2\in\mathcal U( M)$ such that $w_2L(\Th)r_2w_2^*\subseteq L(\Gamma_2)$. Since $L(\Sigma)$ and $L(\Gamma_2)$ are II$_1$ factors, we may assume that $w_2r_2w_2^*\in L(\Sigma)$. 
Let $r_1=1-r_2$.  We will prove that $\Th,r_1,r_2$ satisfy the conclusion of the theorem.

Our first goal is to prove the following:

\begin{claim}\label{conjugate}
There is $u\in\mathcal U( M)$ such that $uL(\Th)r_iu^*\subseteq L(\Gamma_i)$, for every $i\in\overline{1,2}$.
\end{claim}

{\it Proof of Claim \ref{conjugate}.}
To prove the claim, it suffices to show
 that $L(\Theta)r_1$ can be unitarily conjugated into $L(\Gamma_1)$. 
 Indeed, then we can find $w_1\in\mathcal U( M)$ such that $w_1L(\Th)r_1w_1^*\subseteq L(\G_1)$. 
 Since $L(\Gamma_1)$ is a II$_1$ factor and $\uptau(r_1)=\uptau(1-w_2r_2w_2^*)$, we may moreover assume that $w_1r_1w_1^*=1-w_2r_2w_2^*$. Since 
$w_2L(\Th)r_2w_2^*\subseteq L(\Gamma_2)$, it is now clear that $u=w_1r_1+w_2r_2$ is a unitary operator which satisfies the claim. 

Towards showing that $L(\Theta)r_1$ can be unitarily conjugated into $L(\Gamma_1)$, let $q\in \mathcal Z(L(\Theta))r_1$ be a non-zero projection. 
As $L(\Delta)'\cap M\subseteq L(\Omega)\subseteq L(\Th)$ and $L(\Delta)\subseteq L(\Th)$, we have that
 $L(\Theta)'\cap M\subseteq\mathcal Z(L(\Delta)'\cap M)$. Thus, $q\in \mathcal Z(L(\Delta)'\cap  M)r_1$. 
Since $C_{\La}(\Delta)$ is non-amenable,  $L(\Delta)'\cap M$ has no amenable direct summand. 
Since
$\G$ is bi-exact relative to  $\{\Gamma_1,\Gamma_2\}$,  \cite[Theorem 15.1.5]{BO08} implies that $L(\Delta)q\prec_{ M}L(\Gamma_j)$, for some $j\in\overline{1,2}$.
Since $\Delta$ is non-amenable and $\Sigma$ is amenable, we have that $L(\Delta)\nprec_{ M}L(\Sigma)$. By proceeding as in the proof of \cite[Theorem 5.1]{IPP05} it follows that
we can find a non-zero projection 
$r\in \mathcal Z(L(\Delta)'\cap  M)q\subseteq L(\Th)q$ and  $v\in\mathcal U( M)$ such that $vL(\Delta)rv^*\subseteq L(\Gamma_j)$.

Since
$L(\Delta\Omega)\subseteq q\mathcal N_{ M}(L(\Delta))''$ and $L(\Th)= \text{W}^*(q\mathcal N^{(1)}_{ M}(L(\Delta\Omega)))$, by using \cite[Lemma 3.5]{Po03} and Lemma \ref{quasi}, we get
$rL(\Delta\Omega)r\subset q\mathcal N_{r M r}(L(\Delta)r)''$ and $rL(\Theta)r\subseteq$ W$^*(q\mathcal N^{(1)}_{r M r}(rL(\Delta\Omega)r))$.
Since $L(\Delta)\nprec_{ M}L(\Sigma)$, by applying Corollary \ref{mixing3} twice, we get that $vrL(\Theta)rv^*\subseteq L(\Gamma_j)$. Since $\Gamma_j$ is icc, this implies that $L(\Theta)z$ can be unitarily conjugated into $L(\Gamma_j)$, where $z\in\mathcal Z(L(\Theta))$ denotes the central support of $r$.
Then we must have that $j=1$. Otherwise, we would get that $z\leqslant r_2$, hence $r\leqslant r_2$, which  contradicts that $r\leqslant q\leqslant r_1$. 
Therefore, $q'=zq$ is a non-zero projection belonging to $\mathcal Z(L(\Th))q$ such that $L(\Theta)q'$ can be unitarily conjugated into $L(\Gamma_1)$. 
Since this statement holds for every non-zero projection $q\in\mathcal Z(L(\Th))r_1$, and $L(\Gamma_1)$ is a II$_1$ factor, we deduce the existence of $u\in\mathcal U( M)$ such that $uL(\Theta)r_1u^*\subseteq L(\Gamma_1)$. 
\hfill$\square$

We continue with the following:
\begin{claim}\label{findex} There is a non-zero projection $r\in (L(\Delta)'\cap  M)r_1$ such that if $e=uru^*\in L(\Gamma_1)$, then  $ ur(L(\Delta) \vee (L(\Delta)'\cap  M) )ru^*\subseteq eL(\Gamma_1)e$ is a finite index inclusion of II$_1$ factors. 
\end{claim}

 \emph{Proof of Claim \ref{findex}}. 
First, let us show that $L(\Gamma_1^1)\prec_{ M}L(\Delta)r_1$.  Otherwise, since $L(\Gamma_1^1)\prec_{M}L(\Delta)$, we would get that $L(\Gamma_1^1)\prec_{ M}L(\Delta)r_2$,  which would imply that $L(\Gamma_1^1)\prec_{ M} L(\Gamma_2)$. 
Then Lemma \ref{intertwiningroups}  would provide $g\in \G$ such that $[\Gamma^1_1: \Gamma_1^1\cap g\Gamma_2 g^{-1}]<\infty$. This contradicts the fact that $\Gamma^1_1$ is non-amenable and $\Gamma_1^1\cap g\Gamma_2g^{-1}<\Sigma$ is amenable, for every $g\in\Gamma$.

Denote $e_1=ur_1u^*\in L(\Gamma_1)$ and $ P:=uL(\Delta)r_1u^*\subseteq e_1L(\Gamma_1)e_1$.
By the previous paragraph we have $L(\Gamma_1^1)\prec_{ M} P$.
Since $\Gamma_1^1$ is non-amenable and $\Sigma$ is amenable, we have $L(\Gamma_1^1)\nprec_{M}L(\Sigma)$. Thus, applying \cite[Theorem 1.1]{IPP05} gives that $L(\Gamma_1^1)\prec_{L(\Gamma_1)} P$.
By \cite[Lemma 2.4(4)]{DHI16} we get that there is a non-zero projection $e_2\in\mathcal Z( P'\cap e_1L(\Gamma_1)e_1)$ such that $L(\Gamma_1^1)\prec_{L(\Gamma_1)} P f$, for any 
non-zero projection $f\in  P'\cap e_1L(\Gamma_1)e_1$ with $f\leq e_2$. Since $L(\Gamma_1^2)=L(\Gamma_1^1)'\cap L(\Gamma_1)$, by applying \cite[Lemma 3.5]{Va07} it follows that $( P'\cap e_1L(\Gamma_1)e_1)e_2\prec^s_{L(\Gamma_1)}L(\Gamma_1^2)$.

Next, let us show that $ P e_2\prec_{L(\Gamma_1)}^sL(\Gamma_1^1)$. Otherwise, we can find a non-zero projection $f\in\mathcal Z( P'\cap e_1L(\Gamma_1)e_1)$ with $f\leqslant e_2$ such that $ P f\nprec_{L(\Gamma_1)}L(\Gamma_1^1)$. On the other hand, the commutant of $ P f$ in $fL(\Gamma_1)f$ contains $uL(C_{\Lambda}(\Delta))r_1u^*$, and thus has no amenable direct summand. By applying \cite[Theorem 15.1.5]{BO08}, we derive that $ P f\prec^s_{L(\Gamma_1)}L(\Gamma_1^2)$. Since $L(\Gamma_1^1)\prec_{L(\Gamma_1)} P f$, by using \cite[Lemma 3.7]{Va07}, we  get that $L(\Gamma_1^1)\prec_{L(\Gamma_1)}L(\Gamma_1^2)$, which is false.

Since $ P e_2\prec_{L(\Gamma_1)}^sL(\Gamma_1^1)$ and $( P'\cap e_1L(\Gamma_1)e_1)e_2\prec^s_{L(\Gamma_1)}L(\Gamma_1^2)$, we get $\mathcal Z( P)e_2\prec^s_{L(\Gamma_1)}L(\Gamma_1^1)$ and $\mathcal Z( P)e_2\prec^s_{L(\Gamma_1)}L(\Gamma_1^2)$. By using \cite[Lemma 2.8(2)]{DHI16}, this implies that
 $\mathcal Z( P)e_2$ is completely atomic.
Using again that $ P e_2\prec_{L(\Gamma_1)}^sL(\Gamma_1^1)$ and \cite[Lemma 3.7]{Va07}, we derive that $L(\Gamma_1^2)\prec_{L(\Gamma_1)}( P'\cap e_1L(\Gamma_1)e_1)f$,
   for any non-zero projection $f\in  P'\cap e_1L(\Gamma_1)e_1$ satisfying $f\leqslant e_2$. In combination with the fact that $(P'\cap e_1L(\Gamma_1)e_1)e_2\prec^s_{L(\Gamma_1)}L(\Gamma_1^2)$,
we similarly get that   $\mathcal Z( P'\cap e_1L(\Gamma_1)e_1)e_2$ is completely atomic.

In conclusion, both $\mathcal Z( P)e_2$ and  $\mathcal Z( P'\cap e_1L(\Gamma_1)e_1)e_2$ are completely atomic. This implies the existence of a non-zero projection $e_3\in  P'\cap e_1L(\Gamma_1)e_1$ with $e_3\leqslant e_2$ such that both $ P e_3$ and $e_3( P'\cap e_1L(\Gamma_1)e_1)e_3$ are II$_1$ factors.
 Since  $ P e_3\prec_{L(\Gamma_1)} L(\Gamma_1^1)$, \cite[Proposition 12]{OP03} then gives a decomposition $e_3L(\Gamma_1)e_3=L(\Gamma_1^1)^{t_1}\bar{\otimes}L(\Gamma_1^2)^{t_2}$, for some $t_1,t_2>0$ with $t_1t_2=\uptau(e_3)$, and a unitary element $w\in e_3L(\Gamma_1)e_3$ such that $w P e_3w^*\subseteq L(\Gamma_1^1)^{t_1}.$

Since $e_3\leqslant e_2$, we get that $L(\Gamma_1^1)^{t_1}\prec_{e_3L(\Gamma_1)e_3} P e_3$. This gives that $L(\Gamma_1^1)^{t_1}\prec_{L(\Gamma_1^1)^{t_1}}w P e_3w^*$.
 From this we get that there is a non-zero projection $e_4\in (w P e_3w^*)'\cap L(\Gamma_1^1)^{t_1}$ such that the inclusion of II$_1$ factors $(w P e_3w^*)e_4\subseteq e_4L(\Gamma_1^1)^{t_1}e_4$ has finite index.  If we put $e=w^*e_4w$, then $e\in  P'\cap e_1L(\Gamma_1)e_1$, $e\leq e_3$, and \begin{equation}\label{index}w P ew^*\subseteq e_4L(\Gamma_1^1)^{t_1}e_4\;\;\;\text{has finite index}.\end{equation}
Since $w(e( P'\cap e_1L(\Gamma_1)e_1)e)w^*=(w P ew^*)'\cap (e_4\otimes 1)(e_3L(\Gamma_1)e_3)(e_4\otimes 1)$ contains $e_4\otimes L(\Gamma_1^2)^{t_2}$,
 we derive that $we( P\vee ( P'\cap e_1L(\Gamma_1)e_1))ew^*$ contains $w P ew^*\bar{\otimes}L(\Gamma_1^2)^{t_2}$. By using \eqref{index}, we get that the inclusion of II$_1$ factors $we( P\vee ( P'\cap e_1L(\Gamma_1)e_1))ew^*\subseteq (e_4\otimes 1)(e_3L(\Gamma_1)e_3)(e_4\otimes 1)$ has finite index. Thus, the inclusion $e( P\vee ( P'\cap e_1L(\Gamma_1)e_1))e\subseteq eL(\Gamma_1)e$ has finite index, which implies that $r=u^*eu$ satisfies  the claim.  $\hfill\square$

Since $L(\Delta)\vee(L(\Delta)'\cap M)\subseteq L(\Delta\Omega)$ we get that $L(\Delta\Omega)'\cap M\subseteq\mathcal Z(L(\Delta)\vee (L(\Delta)'\cap M))$. Thus, Claim \ref{findex} implies that  $rL(\Delta\Omega)r$ is a II$_1$ factor
and the inclusion
$urL(\Delta\Omega)ru^*\subseteq eL(\Gamma_1)e$ has finite index. 
 By using \cite[Proposition 1.3]{PP86} this entails that $$eL(\Gamma_1)e\subseteq  \text{W}^*(q\mathcal{N}^{(1)}_{ur M ru^*}(urL(\Delta\Omega) ru^*)).$$
  Since $urL(\Delta\Omega)ru^*\subseteq eL(\Gamma_1)e$ has no amenable direct summand and $\Sigma$ is amenable, by Corollary \ref{mixing3}(1) we also get the reverse inclusion.
In combination with Lemma \ref{quasi} and equation \eqref{theta}, we conclude that  $$eL(\Gamma_1)e=\text{W}^*(q\mathcal{N}^{(1)}_{ur M ru^*}(urL(\Delta\Omega) ru^*))=ur\;\text{W}^*(q\mathcal{N}^{(1)}_ M(L(\Delta\Omega)))\;ru^*=urL(\Theta)ru^*.$$

 This implies in particular that  $rL(\Delta\Omega)r\subseteq rL(\Theta)r$ is a finite index inclusion of II$_1$ factors. By \cite[Proposition 1.3]{PP86} it follows that $L(\Theta)\prec_{L(\Th)}L(\Delta\Omega)$. By Lemma \ref{intertwiningroups} we get that  $\Delta\Omega<\Theta$ has finite index. In particular, $\text{Comm}^{(1)}_\La(\Delta\Omega)=\text{Comm}^{(1)}_\La(\Theta)$ and since $\text{Comm}^{(1)}_\La(\Delta\Omega)\subseteq \Theta$ we must have that $\text{Comm}^{(1)}_\La(\Theta)=\Theta$.  By Corollary \ref{FGS} we thus have that $q\mathcal N_{M}^{(1)}(L(\Th))=L(\Th)$.
Since $uL(\Th)r_1u^*\subseteq L(\G_1)$, $urL(\Theta)ru^*=eL(\Gamma_1)e$, and $L(\G_1)$ is a II$_1$ factor, by Lemma \ref{cornertowhole} we deduce that $uL(\Th)r_1u^*=ur_1u^*L(\Gamma_1)ur_1u^*$. This finishes the proof of the theorem.
\hfill$\blacksquare$

\noindent
\subsection*{Proof of Theorem \ref{peripheralidentification2}.}
 By combining Theorems \ref{comm1} and \ref{peripheralidentification1}, for every $i\in\overline{1,2}$, we can find a subgroup $\Th_i<\La$, a unitary element $u_i\in  M$, and a non-zero projection $r_i\in \mathcal Z(L(\Th_i))$ such that  $\text{Comm}_{\La}(\Th_i)=\Th_i$, $q_i:=u_ir_iu_i^*\in L(\G_i)$, and
\begin{equation}\label{120}
\begin{array}{ccc}   u_1L(\Th_1)r_1u_1^*= q_1L(\Gamma_1)q_1 && u_1L(\Th_1)(1-r_1)u_1^*\subseteq L(\Gamma_2)\\  u_2L(\Th_2)(1-r_2)u_2^*\subseteq L(\Gamma_1)
&&u_2L(\Th_2)r_2u_2^*=q_2 L(\Gamma_2)q_2.
\end{array}
\end{equation}

If $r_1=r_2=1$, then conclusion (1)  holds. Therefore, in order to complete the proof, it suffices to prove that if either  $r_1\not=1$ or $r_2\not=1$, then conclusion (2) holds. Due to symmetry, we can further reduce to the case when $r_1\not=1$.

Since $r_1\not=1$, we can find a non-zero projection $r\in L(\Th_1)(1-r_1)$ such that $\uptau(r)\leq \uptau(r_2)$. Since $L(\Gamma_2)$ is a II$_1$ factor, \eqref{120} implies that we can find $v\in\mathcal U( M)$ such that $vL(\Th_1)rv^*\subseteq L(\Th_2)r_2$.
Thus, $L(\Th_1)\prec_{ M}L(\Th_2)$. By applying Lemma \ref{intertwiningroups}, we deduce the existence of $h\in\La$ such that  $[\Th_1:\Th_1\cap h\Th_2h^{-1}]<\infty$. Therefore, after replacing $\Th_2$ with $h\Th_2h^{-1}$, we may assume that in addition to \eqref{120} we also have that $[\Th_1:\Th]<\infty$, where $\Th:=\Th_1\cap\Th_2$. 
In particular, since $\Th_1$ is non-amenable, $\Th$ is non-amenable.

Next, we claim that $r_1r_2=(1-r_1)(1-r_2)=0$. Otherwise, by using \eqref{120} and applying Lemma \ref{intersect}, it follows that we can find $g\in\Gamma$ such that $L(\Th)\prec_{ M}L(\Gamma_1\cap g\Gamma_2g^{-1})$. Since $\Gamma_1\cap g\Gamma_2g^{-1}\subseteq\Sigma$, $\Sigma$ is amenable, and $\Th$ is non-amenable, this leads to a contradiction.

Now, the claim implies that $r_1+r_2=1$. Thus, by \eqref{120} we have that $u_1L(\Th)r_1u_1^*\subseteq L(\Gamma_1)$ and $u_2L(\Th)r_1u_2^*\subseteq L(\Gamma_1)$. Since $\Th$ is non-amenable and $\Sigma$ is amenable,  $L(\Th)\nprec_{ M}L(\Sigma)$. By Lemma \ref{mixing3}(1), we get that $u_1r_1u_2^*\in L(\Gamma_1)$.  In combination with \eqref{120}, this implies that $u_1L(\Theta_2)r_1u_1^*\subseteq L(\Gamma_1)$, hence $u_1L(\Th_2)r_1u_1^*\subseteq q_1L(\G_1)q_1$. Similarly, we get that $u_1r_2u_2^*\in L(\G_2)$. Hence, if we put $\tilde q_2:=u_1r_2u_1^*$, then $\tilde q_2\in L(\Gamma_2)$ and \eqref{120} gives that $u_1L(\Th_2)r_2u_1^*=\tilde q_2L(\G_2)\tilde q_2$.

Since $u_1L(\Th_2)r_1u_1^*\subseteq q_1L(\G_1)q_1$ and $u_1L(\Th_1)r_1u_1^*=q_1L(\G_1)q_1$, we get that $L(\Th_2)r_1\subseteq L(\Th_1)r_1$. This implies that $v_hr_1=0$, for all $h\in \Th_2\setminus\Th_1$. Since $v_hr_1\not=0$, for every $h\in\La$, we conclude that $\Th_2\subseteq\Th_1$ and so $\Th=\Th_2$.  In particular, the inclusion $\Th_2<\Th_1$ has finite index and therefore $\Th_1=\text{Comm}_{\La}^{(1)}(\Th_1)=\text{Comm}_{\La}^{(1)}(\Th_2)=\Th_2$. It follows that  $\Th=\Th_1=\Th_2$ satisfies (2). 
\hfill$\blacksquare$

In the next section, we will prove that if $\Sigma$ is icc and has trivial one-sided commensurator in $\G_1$ and $\G_2$, then condition (2) from Theorem \ref{peripheralidentification2} can be ruled out (see Proposition \ref{noniso}). Here, we point out another general situation in which this is the case.

\begin{cor}\label{BS}
Let  $\G=\G_1\ast_{\Sigma}\G_2$ be as in Assumption \ref{AFP}, and put $ M=L(\G)$. 
Assume  additionally that either $\G_1$ has property (T) and $\G_2$ does not, or that $\G_1$ has Haagerup's property and $\G_2$ does not.
Let $\La$ be an arbitrary group such that $ M=L(\Lambda)$. 

Then for any $i\in\overline{1,2}$, there exists a subgroup $\La_i<\La$ such that $u_i L(\La_i)u^*_i=L(\Gamma_i)$ for $u_i\in\mathcal U( M)$.

\end{cor}

{\it Proof.}
Using the assumptions made on $\G_1$ and $\G_2$, Proposition \ref{HT} guarantees that condition (2) from Theorem \ref{peripheralidentification2} does not hold. The conclusion now follows from Theorem \ref{peripheralidentification2} .
\hfill$\blacksquare$

\section{Proof of Theorem \ref{A}}\label{4}
This section is devoted to the proof of Theorem \ref{A}, whose setup we now recall. 
Let $ M=L(\Gamma)$, where $\G=\Gamma_1\ast_{\Sigma}\Gamma_2$ is an amalgamated free product group satisfying the following conditions:

\begin{enumerate} 
\item $\Sigma$ is an icc amenable group and $\text{Comm}_{\G_i}^{(1)}(\Sigma)=\Sigma$, for every  $i\in\overline{1,2}$.
\item $\Gamma_i=\Gamma_i^1\times\Gamma_i^2$, where $\Gamma_i^j$ is an icc, non-amenable, bi-exact group, for every $i,j\in\overline{1,2}$. 
\end{enumerate}

 In order to derive Theorem \ref{A}, we will need the following result, whose proof we postpone until the end of this section.

\begin{prop}\label{noniso}
Let $\La$ be a countable group such that $ M=L(\La)$.

Then there do not exist  a subgroup $\Th<\La$, non-zero projections $r_1,r_2\in \mathcal Z(L(\Th))$, and a unitary $u\in M$, such that $r_1+r_2=1$, $q_i:=ur_iu^*\in L(\Gamma_i)$ and $uL(\Th)r_iu^*=q_iL(\G_i)q_i$, for every $i\in\overline{1,2}$.
\end{prop}

{\it Proof of Theorem \ref{A}.} Let $\La$ be a group such that $M=L(\La)$.
Denote by $\set{u_g}_{g\in \Gamma}$ and $\set{v_h}_{h\in \Lambda} $ the canonical unitaries generating $ M$.  
Theorem \ref{peripheralidentification2} implies that either condition (1) or (2) from its conclusion hold. By Proposition \ref{noniso}, condition (2) cannot hold. Thus, we deduce that
  for every $i\in\overline{1,2}$, we can find a subgroup $\Th_i<\La$ and $v_i\in\mathcal U( M)$ such that  $v_iL(\Th_i)v_i^*=L(\Gamma_i).$

In particular, we get that $L(\Sigma)= L(\Gamma_1)\cap L(\Gamma_2)=v_1L(\Th_1)v_1^*\cap v_2L(\Th_2)v_2^*$. Lemma \ref{intersect} implies that $L(\Sigma)\prec_{ M}L(\Th_1\cap h\Th_2h^{-1})$, for some $h\in\La$.
Define $\La_1=\Th_1$, $\La_2=h\Th_2h^{-1}$, and $\Delta=\La_1\cap\La_2$. 
Letting $u_1=v_1$ and $u_2=v_2u_h^*$, we have that
\begin{equation}\label{equal} u_iL(\La_i)u_i^*=L(\Gamma_i),\;\;\;\text{for every $i\in\overline{1,2}$.}
\end{equation}
\noindent
In particular, we get that $L(\Delta)\subseteq u_1^*L(\Gamma_1)u_1\cap u_2^*L(\Gamma_2)u_2$. Since $\Gamma_1\cap g\Gamma_2g^{-1}\subseteq \Sigma$, for all $g\in\Gamma$, by applying Lemma \ref{intersect} we conclude that \begin{equation}\label{strongly}L(\Delta)\prec_{ M}^sL(\Sigma).\end{equation}

Since $L(\Sigma)\prec_{ M}L(\Delta)$, by \cite[Lemma 2.4(4)]{DHI16}, there is a non-zero projection $z\in\mathcal Z(L(\Delta)'\cap M)$ such that $L(\Sigma)\prec_{ M}L(\Delta)q'$, for any non-zero projection $q'\in (L(\Delta)'\cap M)z$. We may moreover assume that $z$ is the largest projection belonging to $\mathcal Z(L(\Delta)'\cap M)$ with this property. 

We claim that for every $i\in\overline{1,2}$, $g\in\Gamma\setminus\Sigma$, and $h\in\Gamma\setminus\Gamma_i$  we have that \begin{equation}\label{noint} L(\Delta)z\nprec_{\paM}L(\Sigma\cap g\Sigma g^{-1}) \;\;\text{and}\;\; L(\Delta)z\nprec_{\paM}L(\Gamma_i\cap h\Gamma_i h^{-1}).\end{equation}
 Indeed, if $L(\Delta)z\prec_{M}L(\Sigma\cap g\Sigma g^{-1})$ (respectively, if $L(\Delta)z)\prec_{M}L(\Gamma_i\cap g\Gamma_i g^{-1})$), for some $g\in\Gamma$, then by \cite[Lemma 2.4 (3)]{DHI16} we can find a non-zero projection $q'\in (L(\Delta)'\cap M)z$ such that $L(\Delta)q'\prec_{M}^{s}L(\Sigma\cap g\Sigma g^{-1})$ (resp., $L(\Delta)q'\prec_{M}^{s}L(\Gamma_i\cap g\Gamma_i g^{-1})$). On the other, since $q'\leqslant z$, we have that $L(\Sigma)\prec_{ M}L(\Delta)q'$. By \cite[Lemma 3.7]{Va07} we get that $L(\Sigma)\prec_{ M}L(\Sigma\cap g\Sigma g^{-1})$ (resp., $L(\Sigma)\prec_ML(\Gamma_i\cap g\Gamma_i g^{-1})$). Lemma \ref{mixing4}(2) gives that $g\in\Sigma$ (resp., $g\in\Gamma_i$), which proves  \eqref{noint}.

\begin{claim}\label{core}
We can find $u\in\mathcal U( M)$ such that $uL(\Delta)zu^*\subseteq L(\Sigma)$.
\end{claim}

\noindent
{\it Proof of Claim \ref{core}.} Let $q'\in (L(\Delta)'\cap M)z$ be a non-zero projection. Since $L(\Delta)q'\prec_{ M}L(\Sigma)$, we can find projections $q\in L(\Delta)$, $r\in L(\Sigma)$, a non-zero isometry $w\in r M qq'$ and a $*$-homomorphism $\varphi:qL(\Delta)qq'\rightarrow rL(\Sigma)r$ satisfying $\varphi(x)w=wx$, for all $x\in qL(\Delta)qq'$.  Moreover, we may assume that $r$ is equal to the support projection of $E_{L(\Sigma)}(ww^*)$.
Put $ P:=\varphi(qL(\Delta)qq')\subseteq rL(\Sigma)r$.

Note that $ P\nprec_{L(\Sigma)}L(\Sigma\cap g\Sigma g^{-1})$, for all $g\in\Gamma\setminus\Sigma$.
Otherwise, \cite[Lemma 1.12]{IPP05} would imply that $qL(\Delta)qq'\prec_{ M}L(\Sigma\cap g\Sigma g^{-1})$, which contradicts \eqref{noint}.
By applying Lemma \ref{mixing} we deduce that   $P'\cap rL(\Sigma)r\subseteq L(\Sigma)$, hence $ww^*\in L(\Sigma)$. Since $w^*w\in q'(L(\Delta)'\cap M)q'q$, we can find a projection $q_0\in q'(L(\Delta)'\cap M)q'$ such that $w^*w=qq_0$.  Thus, $w(qL(\Delta)qq_0)w^*\subseteq L(\Sigma)$.  Let $z_0$ be the central support of $q$ in $L(\Delta)$. Since $L(\Sigma)$ is a II$_1$ factor, it follows  that there is $\eta\in \mathcal U( M)$ such that $\eta L(\Delta)z_0q_0\eta^*\subseteq L(\Sigma)$.

Thus, for any non-zero projection $q'\in (L(\Delta)'\cap M)z$, we found a non-zero projection $q''$ in $q'(L(\Delta)'\cap M)q'$ such that $L(\Delta)q''$ can be unitarily conjugated into $L(\Sigma)$. Since $L(\Sigma)$ is a II$_1$ factor, the claim follows by a maximality argument 
 (see the proof of \cite[Theorem 5.1]{IPP05}). \hfill$\square$

Next, we define by $\Omega<\La$ the subgroup generated by $\text{Comm}_{\La}^{(1)}(\Delta)$, and prove the following:

\begin{claim}\label{core2}We have that $uzL(\Omega)zu^*\subseteq L(\Sigma)$ and there is a non-zero projection $z'\in zL(\Om)z$ such that $uz'L(\Om)z'u^*=pL(\Sigma)p$, where $p=uz'u^*$.
\end{claim}

{\it Proof of Claim \ref{core2}.}
Put $e:=uzu^*\in L(\Sigma)$. First, since $L(\Delta)z\nprec L(\Sigma\cap g\Sigma g^{-1})$, for every $g\in\Gamma\setminus\Sigma$,
by Lemma \ref{mixing} we deduce that $q\mathcal N_{e M e}^{(1)}(uL(\Delta)zu^*)\subseteq eL(\Sigma)e$. On the other hand, the combination of Lemma \ref{quasi} and Corollary \ref{FGS}
yields that $uzL(\Omega)zu^*=W^*(q\mathcal N_{e M e}^{(1)}(uL(\Delta)zu^*))$. Putting together the last two facts, we deduce that indeed $uzL(\Omega)zu^*\subseteq L(\Sigma)$.

Second, put $ Q:=uL(\Delta)zu^*\subseteq eL(\Sigma)e$.  Since $L(\Sigma)$ is a II$_1$ factor and $L(\Sigma)\prec_{M}L(\Delta)z$, we get that $eL(\Sigma)e\prec_{ M} Q $. Thus, we can find projections $p\in eL(\Sigma)e$, $q\in Q$, a non-zero partial isometry $v\in q M p$, and a $*$-homomorphism $\theta:pL(\Sigma)p\rightarrow q Q q$ such that $\theta(x)v=vx$, for all $x\in pL(\Sigma)p$. Since $q Q q\subseteq L(\Sigma)$, it follows that $v\in W^*(q\mathcal N^{(1)}_{ M}(L(\Sigma)))$. Since $\text{Comm}_{\G_i}^{(1)}(\Sigma)=\Sigma$, for all $i\in\overline{1,2}$, combining Corollary \ref{FGS} and Lemma \ref{mixing4}(1) yields that $v\in L(\Sigma)$. Thus, after shrinking $p$, we may assume that $v^*v=p$. 
Moreover, since $L(\Sigma)$ is a II$_1$ factor and $ Q \subseteq eL(\Sigma)e$ is diffuse, we may assume that $p\in Q$.  

Now, if $x\in pL(\Sigma)p$, then $vx(p Q p)\subseteq vpL(\Sigma)p\subseteq (qQq)v$. This implies that $vx\in W^*(q\mathcal N_{eL(\Sigma)e}^{(1)}(Q))$ (see the proofs of \cite[Lemma 3.5]{Po03} or Lemma \ref{cornertowhole}).
Thus,  $pL(\Sigma)p\subseteq W^*(q\mathcal N_{eM)e}(Q))$. 
On the other hand,  the moreover assertion of Lemma \ref{cornertowhole} and Lemma \ref{mixing2} imply that $$q\mathcal N_{e M e}^{(1)}( Q)=uz\; q\mathcal N_{ M}^{(1)}(L(\Delta)) \;zu^*\subseteq uzL(\Om)zu^*.$$ By combining the last two
 inclusions we deduce that $pL(\Sigma)p\subseteq p (uzL(\Om)zu^*)p$. Therefore, if $z'\in L(\Delta)z\subseteq zL(\Om)z$ is such that $p=uz'u^*$, then $pL(\Sigma)p\subseteq uz'L(\Om)z'u^*$. Since the reverse inclusion also holds, the second assertion of the claim follows.
\hfill$\square$

Before finishing the proof, we need one final claim:

\begin{claim}\label{core3} $[\Om:\Delta]<\infty$ and there is $w\in\mathcal U( M)$  such that $wL(\Om)w^*=L(\Sigma)$.
\end{claim} 

{\it Proof of Claim \ref{core3}.} 
Note that $vuzL(\Omega)zu^*\subseteq vpL(\Sigma)p=\theta(pL(\Sigma)p)v\subseteq uL(\Delta)zu^*v$. Thus, if $\xi=u^*vu$, then $\xi\in z M z$ and $\xi zL(\Omega)z\subseteq L(\Delta)z\xi$. In particular, $\xi L(\Delta)\subseteq L(\Delta)\xi$. Corollary \ref{FGS} implies that $\xi\in zL(\Om)z$. Therefore, $L(\Om)\prec_{L(\Om)}L(\Delta)$, and Lemma \ref{intertwiningroups} gives that $[\Om:\Delta]<\infty$.

Since $[\Om:\Delta]<\infty$, 
 $\text{Comm}^{(1)}_{\La}(\Omega)=\text{Comm}^{(1)}_{\La}(\Delta)\subseteq\Omega$ and thus $\text{Comm}^{(1)}_{\La}(\Omega)=\Omega$.
Corollary \ref{FGS} gives that $q\mathcal N_{ M}^{(1)}(L(\Omega))=L(\Omega)$.
 Moreover, since $L(\Delta)\prec_{ M}^sL(\Sigma)$ by \eqref{strongly}, by using again that $[\Omega:\Delta]<\infty$, we get that $L(\Omega)\prec_{ M}^{s}L(\Sigma)$. Since  $L(\Sigma)$ is a II$_1$ factor and $uz'L(\Omega)z'u^*=pL(\Sigma)p$, we can apply Lemma \ref{cornertowhole} and deduce the existence of $w\in\mathcal U( M)$ such that $wL(\Om)w^*=L(\Sigma)$.
\hfill$\square$
 
We are now ready to finish the proof. Let $q'\in L(\Delta)'\cap M$ be a non-zero projection. Then $q'\in L(\Omega)$ and since $[\Omega:\Delta]<\infty$ we have that $q'L(\Omega)q'\prec_{ M}^sL(\Delta)q'$. Moreover,  $L(\Sigma)\prec_{ M}q'L(\Omega)q'$ by Claim \ref{core2}, hence \cite[Lemma 3.7]{Va07} allows us to conclude that $L(\Sigma)\prec_{ M}L(\Delta)q'$. Since this holds for any non-zero projection $q'\in L(\Delta)'\cap M$,  the maximality 
 property of $z$  implies that $z=1$. 

By Claim \ref{core2} we get that  $ Q:=uL(\Delta)u^*\subseteq L(\Sigma)$. Let $i\in\overline{1,2}$. Since $z=1$, \eqref{noint} implies that  $L(\Delta)\nprec_{ M}L(\Gamma_i\cap g\Sigma g^{-1})$, for all $g\in\Gamma\setminus\Gamma_i$.
Since $u_iu^* Q=u_iL(\Delta)u^*\subseteq u_iL(\La_i)u^*=L(\G_i)u_iu^*$
 by equation \eqref{equal},  Lemma \ref{mixing} gives that $u_iu^*\in L(\Gamma_i)$. Thus, we get that $$u^*L(\Gamma_i)u=u^*(uu_i^*)L(\Gamma_i)(uu_i^*)^*u=u_i^*L(\Gamma_i)u_i=L(\Lambda_i).$$ 
Therefore, $L(\Delta)=L(\La_1)\cap L(\La_2)=u^*(L(\G_1)\cap L(\G_2))u=u^*L(\Sigma)u$. This finishes the proof.
\hfill$\blacksquare$

\subsection*{Proof of Proposition \ref{noniso}.}
  Assume by contradiction that the conclusion of Proposition \ref{noniso} is false. 
After replacing $\La$ with $u\La u^*$, we find a group $\La$ satisfying $ M=L(\La)$, a subgroup $\Th<\La$  and non-zero projections $r_1,r_2\in \mathcal Z(L(\Th))\cap L(\Sigma)$ such that  $r_1+r_2=1$ and \begin{equation}\label{ri}L(\Th)r_i=r_iL(\G_i)r_i,\;\;\;\text{for all $i\in\overline{1,2}$}.\end{equation}

Since $L(\Sigma)$ is a II$_1$ factor, there is a non-zero partial isometry $v\in L(\Sigma)$ such that $vv^*\leqslant r_1$ and $v^*v\leqslant r_2$. Then  $vv^*L(\Sigma)vv^*\subseteq r_1L(\Sigma)r_1\cap vr_2L(\Sigma)r_2v^*\subseteq L(\Th)\cap vL(\Th)v^*$. 
The moreover assertion of Lemma \ref{intersect} implies that there exists $h\in\La$ such that $L(\Sigma)\prec_{ M}L(\Th\cap h\Th h^{-1})$ and $\uptau(vv_h^*)\not=0$. Moreover, since $v=r_1vr_2$, we get that $E_{L(\Th)}(v)=r_1E_{L(\Th)}(v)r_2=r_1r_2E_L{(\Th)}(v)=0$, hence $h\in\La\setminus\Th$.  
Thus, if we put $\Delta:=\Th\cap h\Th h^{-1}$, then \begin{equation}\label{delta1}L(\Sigma)\prec_{ M}L(\Delta).\end{equation}
We claim that \begin{equation}\label{delta2} L(\Delta)\prec_{ M}^{s}L(\Sigma). \end{equation} For this, let $p\in\mathcal Z( (L(\Delta)r_1)'\cap r_1L(\Gamma_1)r_1)$ be the largest projection such that $L(\Delta)p\nprec_{L(\G_1)}L(\Sigma)$. 
First, since we have $L(\Delta)p(v_hr_1)=pL(\Delta)(v_hr_1)=pv_hL(h^{-1}\Delta h)r_1\subseteq pv_h L(\Th)r_1\subseteq pv_hL(\G_1)$, Corollary \ref{mixing3} allows us to conclude that $pv_hr_1\in r_1L(\G_1)r_1$. In particular, $pv_hr_1\in L(\Th)$ and since $r_1, p\in L(\Th)$ while $h\in\La\setminus\Th$, we get that $pv_hr_1=pE_{L(\Th)}(v_h)r_1=0$.
Secondly, since $L(\Delta)p(v_hr_2)=p L(\Delta)(v_hr_2)= pv_hL(h^{-1}\Delta h)r_2\subseteq pv_hL(\Th)r_2\subseteq pv_hL(\Gamma_2)$ and we have $L(\Delta)p\nprec_{L(\G_1)}L(\Sigma)$, \cite[Theorem 1.1]{IPP05} implies that  $pv_hr_2=0$.

Since $r_1+r_2=1$,   the last paragraph gives that $p=0$.  This implies that $L(\Delta)r_1\prec^{s}_{L(\G_1)}L(\Sigma)$. 
 Similarly, it follows that  $L(\Delta)r_2\prec_{L(\G_2)}^{s}L(\Sigma)$. These together prove \eqref{delta2}.
 
 Let $\Omega<\La$ be the subgroup generated by $\text{Comm}_{\La}^{(1)}(\Delta)$.
In the proof of Theorem \ref{A}, we showed that if $\Delta<\La$ satisfies conditions \eqref{delta1} and \eqref{delta2}, then $[\Omega:\Delta]<\infty$ and $wL(\Om)w^*=L(\Sigma)$, for some $w\in\mathcal U( M)$ (see  Claim \ref{core3}).
In particular, $\Omega$ is icc.

Put $ Q:=wL(\Delta)w^*\subseteq L(\Sigma)$. Since $\Delta\subseteq\Om$, we have $r_iw^* Q=r_iL(\Delta)w^*\subseteq r_iL(\G_i)r_iw^*$, for all $i\in\overline{1,2}$.
Note that $ Q\nprec_{ M}L(\Gamma_i\cap g\Gamma_i g^{-1})$, for all $g\in\Gamma_i\setminus\Sigma$. Otherwise, since $[\Omega:\Delta]<\infty$ we have that $L(\Sigma)\prec_{ M} Q 	q$, for any non-zero projection $q\in Q'\cap  M$, and  \cite[Lemma 3.7]{Va07} would imply that $L(\Sigma)\prec_{ M}L(\Gamma_i\cap g\Gamma_i g^{-1})$. This however contradicts Lemma \ref{mixing4}(3). 

We can therefore apply Lemma \ref{mixing} to derive that $r_iw^*\in L(\G_i)$.
Let $p_i=wr_iw^*\in\mathcal P(L(\G_i))$.
Then $p_1,p_2$ are non-zero and since $p_1+p_2=1$ we get that $p_1,p_2\in L(\Sigma)$. Moreover, we have that \begin{align*}wL(\Th)w^*=w(L(\Th)r_1\oplus L(\Th)r_2)w^*&=wr_1L(\G_1)r_1w^*\oplus wr_2L(\G_2)r_2w^*\\ &=p_1L(\G_1)p_1\oplus p_2L(\G_2)p_2.\end{align*}
From this we deduce that \begin{align*}wL(\Om\cap\Th)w^*=wL(\Om)w^*\cap wL(\Th)w^*&=L(\Sigma)\cap (p_1L(\G_1)p_1\oplus p_2L(\G_2)p_2)\\&=p_1L(\Sigma)p_1\oplus p_2L(\Sigma)p_2.\end{align*}
In particular, since $p_1,p_2$ are non-zero, we conclude that $L(\Om\cap\Th)$ is not a factor and therefore $\Om\cap\Th$ is not icc. On other hand, since $[\Om:\Delta]<\infty$ and $\Delta\subseteq\Om\cap\Th$, we get that $[\Om:\Om\cap\Th]<\infty$. This altogether contradicts the fact that $\Om$ is icc.
\hfill$\blacksquare$


\section{Proof of Corollary \ref{B}}

In this section, we prove Corollary \ref{B}. Its proof relies on Theorem \ref{A} and the following result:

\begin{theorem}\label{super}
Let $\Gamma_1, \Gamma_2$ be icc, non-amenable, bi-exact groups. Put $\Gamma=\Gamma_1\times\Gamma_2$ and $M=L(\G)$. Let $\Sigma$ be an icc group and $\pi_i:\Sigma\rightarrow\Gamma_i$ an injective homomorpism such that $\{\pi_i(g)h\pi_i(g)^{-1}|g\in\Sigma\}$ is infinite, for all $h\in\Gamma_i\setminus\{e\}$ and $i\in\overline{1,2}$. We identify $\Sigma$ with $\{(\pi_1(g),\pi_2(g))|g\in\Sigma\}<\Gamma$. 

Let $\Delta<\La$ be countable groups such that $ M=L(\La)$ and $L(\Sigma)=L(\Delta)$. 

Then we can find a decomposition $\La=\La_1\times\La_2$ and a unitary $u\in M$ such that
$$ \mathbb T\Sigma=u\mathbb T\Delta u^*\;\;\;\text{and}\;\;\; L(\G_i)=uL(\La_i)u^*,\;\;\;\text{for all $i\in\overline{1,2}$.}$$
\noindent
\end{theorem}

Recall from \cite[Section 4]{Io10}, that the {\it height} of an element $x\in L(\La)$ is defined as $$h_{\La}(x)=\max_{h\in\La}|\uptau(xv_h^*)|.$$
In the proof of Theorem \ref{super}, we will make crucial use of \cite[Theorem 3.1]{IPV10}. This asserts that if $\G$ is any countable group such that $L(\G)=L(\La)$ and $$\inf_{g\in\G}h_{\La}(u_g)>0,$$ then there is a unitary $u\in L(\G)=L(\La)$ such that $\mathbb T\G=u\mathbb T\La u^*.$

\subsection*{Proof of Theorem \ref{super}}
By \cite[Corollary B]{CdSS15} we can find a decomposition $\Lambda=\Lambda_1\times\Lambda_2$, where $\Lambda_1,\Lambda_2$ are icc groups, $t_1,t_2>0$ with $t_1t_2=1$, and $x\in\mathcal U( M)$ such that $L(\La_1)=xL(\G_1)^{t_1}x^*$ and $L(\La_2)=xL(\G_2)^{t_2}x^*$. 
Let $d\geqslant \max\{t_1,t_2\}$ be an integer. For $i\in\overline{1,2}$, let $p_i\in\mathbb M_d(L(\La_i))$ be a projection with $(\uptau\otimes\text{Tr})(p_i)=t_i$, where $\text{Tr}$ denotes the non-normalized trace on $\mathbb M_d(\mathbb C)$. Then the above implies that we can find a partial isometry $v\in \mathbb M_d(L(\La_1))\bar{\otimes}\mathbb M_{d}(L(\La_2))$ such that $vv^*=e_{1,1}\otimes e_{1,1}$, $v^*v=p_1\otimes p_2$, where $e_{1,1}\in\mathbb M_d(\mathbb C)$ denotes the elementary matrix corresponding to the $(1,1)$-entry, and 
if we identify $L(\La_i)\equiv e_{1,1}\mathbb M_d(L(\La_i))e_{1,1}$ in the natural way, then
\begin{equation}\label{identify}
L(\La_1)\otimes e_{1,1}=v(p_1\mathbb M_d(L(\G_1))p_1\otimes p_2)v^*\;\;\;\text{and}\;\;\; e_{1,1}\otimes L(\La_2)=v(p_1\otimes p_2\mathbb M_{d}(L(\G_2))p_2)v^*.
\end{equation}

Let $\rho_i$ be the restriction of the projection $\Lambda\rightarrow\La_i$ to $\Delta$.
We claim that
$\rho_i$ is one-to-one, for all $i\in\overline{1,2}$.
We only treat the case $i=1$, since the case $i=2$ is similar.
To this end, let $\Om=\ker(\rho_1)$. Since $\Delta$ is icc, in order to show that $\Om=\{e\}$, it suffices to prove that $\Om$ is finite. Assume  that $\Om$ is infinite, and 
let $h_n\in\Om$ be a sequence satisfying $h_n\rightarrow\infty$. 
Since $v_{h_n}\in 1\otimes L(\Lambda_2)$, \eqref{identify} implies the existence of $T\subseteq\G_1$ finite such that $\|v_{h_n}-e(v_{h_n})\|_2\leqslant 1/2$, for all $n\geqslant 1$, where $e$ denotes the orthogonal projection from $\ell^2(\Gamma)=\ell^2(\G_1)\otimes\ell^2(\G_2)$ onto the closed linear span of $\{u_{g}\otimes L(\G_2)|g\in T\}$. On the other hand, $v_{h_n}\in L(\Delta)=L(\Sigma)$, for all $n\geqslant 1$. Thus, we get that $$\|e(v_{h_n})\|_2^2=\sum_{g\in\pi_1^{-1}(T)}|\uptau(v_{h_n}u_g^*)|^2,\;\;\;\text{for all $n$}.$$
Since $\pi_1$ is one-to-one and $T$ is finite, $\pi_1^{-1}(T)\subseteq\Delta$ is finite. Since $v_{h_n}\rightarrow 0$, weakly, we conclude that $\|e(v_{h_n})\|_2\rightarrow 0$, as $n\rightarrow\infty$, which gives a contradiction, and proves the claim.

We continue by establishing the following:
\begin{claim}\label{height}
$\inf_{g\in\Sigma}h_{\Delta}(u_g)>0$.
\end{claim}

{\it Proof of Claim \ref{height}.}
Using \eqref{identify}, for every $i\in\overline{1,2}$, we can find a finite set $S_i\subseteq\La_i$ such that for every $u_i\in\mathcal U(L(\Gamma_i))$, there is $v_i$ in the linear span of $\{v_{h_1}\otimes v_{h_2}| h_1\in\La_1,h_2\in\La_2, h_i\in S_i\}$  satisfying $\|v_i\|\leqslant 1$ and $\|u_i-v_i\|_2\leqslant 1/8$.

Let $g\in\Sigma$. Then for every $i\in\overline{1,2}$ we can find $v_i\in (M)_1$ such that $\|u_{\pi_i(g)}-v_i\|_2\leqslant 1/8$ and $$v_1=\sum_{(h_1,h_2)\in\La_1\times S_2}c_{h_1,h_2}(v_{h_1}\otimes v_{h_2})\;\;\;\text{and}\;\;\;v_2=\sum_{(h_1',h_2')\in S_1\times \La_2}d_{h_1',h_2'}(v_{h_1'}\otimes v_{h_2'}),$$ for some $c_{h_1,h_2}, d_{h_1',h_2'}\in\mathbb C$. 

Since $u_g=u_{\pi_1(g)}u_{\pi_2(g)}\in L(\Delta)$, we get that $\|u_g-v_1v_2\|_2\leqslant 1/4$, hence $\|u_g-E_{L(\Delta)}(v_1v_2)\|_2\leqslant 1/4$.
Write $u_g=\sum_{k\in\Delta}a_kv_k$, and notice that $E_{L(\Delta)}(v_1v_2)=\sum_{k\in\Delta}b_kv_k$, where $$b_k=\sum_{\substack{
 (h_1,h_2)\in\La_1\times S_ 2, (h_1',h_2')\in S_1\times\La_2 \\
   h_1h_1'=\rho_1(k),h_2h_2'=\rho_2(k)}} c_{h_1,h_2}d_{h_1',h_2'}.$$

Now, fix $(h_1,h_2)\in\La_1\times S_2$. If $k\in\Delta$ is such that there is $(h_1',h_2')\in S_1\times\La_2$ satisfying $h_1h_1'=\rho_1(k)$, then $k\in\rho_1^{-1}(h_1S_1)$. Since $\rho_1$ is one-to-one, there are at most $|S_1|$ such $k\in\Delta$. Similarly, given $(h_1',h_2')\in S_1\times\La_2$, there are at most $|S_2|$ elements $k\in\Delta$ for which there is $(h_1,h_2)\in\La_1\times S_2$ satisfying $h_2h_2'=\rho_2(k)$.
 Using the inequality $|cd|\leqslant c^2+d^2$, for all $c,d\in\mathbb R$, we conclude that 
\begin{equation}\label{bk}
\sum_{k\in\Delta}|b_k|\leqslant \sum_{(h_1,h_2)\in \La_1\times S_2}|S_1|\;c_{h_1,h_2}^2+\sum_{(h_1',h_2')\in S_1\times\La_2}|S_2|\;d_{h_1',h_2'}^2\leqslant |S_1|+|S_2|.
\end{equation}
Next, let $T=\{k\in\Delta||a_k-b_k|\leqslant |a_k|/2\}$. Then \begin{equation}\label{ak}\sum_{k\in\Delta\setminus T}|a_k|^2\leqslant\sum_{k\in\Delta}4|a_k-b_k|^2=4\|u_g-E_{L(\Delta)}(v_1v_2)\|_2^2\leqslant 1/4.  \end{equation}
On the other hand, if $k\in T$, then $|a_k|\leqslant 2|b_k|$, and thus by using \eqref{bk} we get that
\begin{equation}\label{akbk}\sum_{k\in T}|a_k|\leqslant 2\sum_{k\in\Delta}|b_k|\leqslant 2(|S_1|+|S_2|).
\end{equation}
Finally, by combining \eqref{ak} and \eqref{akbk},  we deduce that 
\begin{align*}1=\sum_{k\in\Delta}|a_k|^2&=\sum_{k\in T}|a_k|^2+\sum_{k\in\Delta\setminus T}|a_k|^2\\&\leqslant h_{\Delta}(u_g)\sum_{k\in T}|a_k|+1/4\\&\leqslant 2(|S_1|+|S_2|)\;h_{\Delta}(u_g)+1/4. \end{align*}
Therefore, we have $h_{\Delta}(u_g)\geqslant (3|S_1|+|S_2|)/8>0$, for any $g\in\Sigma$. This proves the claim.
\hfill$\square$

We are now ready to finish the proof.
First,  combining Claim \eqref{height} and \cite[Theorem 3.1]{IPV10} allows us to deduce the existence of $w\in L(\Sigma)$, an isomorphism $\delta:\Sigma\ra\Delta$, and a character $\eta:\Sigma\ra\mathbb T$ such that $u_g=\eta(g)w v_{\delta(g)} w^*$, for all $g\in\Sigma$. Moreover, after replacing $\La$ with $w\La w^*$, we may assume that $w=1$. In other words, we have
\begin{equation}\label{uv} u_{\pi_1(g)}\otimes u_{\pi_2(g)}=u_g=\eta(g)v_{\delta(g)}=\eta(g)(v_{\rho_1(\delta(g))}\otimes v_{\rho_2(\delta(g))}),\;\;\;\text{for all $g\in\Gamma$.}\end{equation}

By equation \eqref{identify}, for every $i\in\overline{1,2}$, we have a homomorphism $\sigma_i:\Sigma\rightarrow\mathcal U(p_i\mathbb M_d(L(\G_i))p_i)$ such that $v_{\rho_1(\delta(g))}\otimes e_{1,1}=v(\sigma_1(g)\otimes p_2)v^*$ and $e_{1,1}\otimes v_{\rho_2(\delta(g))}=v(p_1\otimes\sigma_2(g))v^*$, for all $g\in\Sigma$.
In combination with \eqref{uv}, we deduce that
\begin{equation}\label{utov2}
u_{\pi_1(g)}\otimes u_{\pi_2(g)}=\eta(g)\; v(\sigma_1(g)\otimes\sigma_2(g))v^*,\;\;\;\text{for all $g\in\Sigma$}.
\end{equation}
For $i\in\overline{1,2}$, we define a unitary representation $\alpha_i:\Sigma\rightarrow\mathcal U(L^2(e_{1,1}\mathbb M_d(L(\G_i))p_i))$ by letting $\alpha_i(g)(\xi)=u_{\pi_i(g)}\xi\sigma_i(g)^*$. Then \eqref{utov2} implies that $(\alpha_1(g)\otimes\alpha_2(g))(v)=\eta(g)v$, for all $g\in\Sigma$.
Therefore, both $\alpha_1$ and $\alpha_2$ are not weakly mixing. 

We continue by using an argument from the proof of  \cite[Lemma 2.5]{PS03}. Let $i\in\overline{1,2}$.
Since $\alpha_i$ is not weakly mixing, we can find an $\alpha_i(\Sigma)$-invariant subspace $\{0\}\not=\mathcal H_i\subseteq L^2(e_{1,1}\mathbb M_d(L(\G_i))p_i)$. Let $\mathcal B_i$ be an orthonormal basis of $\mathcal H_i$. Then  $\xi=\sum_{\zeta\in\mathcal B_i}\zeta\zeta^*\in L^1(e_{1,1}\mathbb M_d(L(\G_i))e_{1,1})=L^1(L(\G_i))$ is independent of the choice of the basis. Since $\{\alpha_i(g)(\zeta)\}_{\zeta\in\mathcal B_i}$ is a basis of $\mathcal H_i$, we get that 
$$\xi=\sum_{\zeta\in\mathcal B_i}\alpha_i(g)(\zeta)\alpha_i(g)(\zeta)^*=u_{\pi_i(g)}\xi u_{\pi_i(g)}^*,\;\;\;\text{for all $g\in\Sigma$.}$$
Since $\{\pi_i(g)h\pi_i(g)^{-1}|g\in\Sigma\}$ is infinite, for all $h\in\G_i\setminus\{e\}$,  this forces $\xi\in\mathbb C1$.  In particular, we derive that $\zeta\in L(\G_i)$, for all $\zeta\in\mathcal B_i$, and thus $\mathcal H_i\subseteq L(\G_i)$. Let $\mathcal K_i\subseteq L(\G_i)$ be the linear span $\{\zeta_1\zeta_2^*|\zeta_1,\zeta_2\in\mathcal H_i\}$. Then $\mathcal K_i$ is  a finite dimensional space which is invariant under the unitary representation $\tau_i:\Sigma\rightarrow\mathcal U(L^2(\G_i))$ given by $\tau_i(g)(\zeta)=u_{\pi_i(g)}\zeta u_{\pi_i(g)}^*$.
Using again the fact that $\{\pi_i(g)h\pi_i(g)^{-1}|g\in\Sigma\}$ is infinite, for all $h\in\G_i\setminus\{e\}$, we deduce that $\mathcal K_i\subseteq\mathbb C1$. 

This implies the existence of a partial isometry  $\omega_i\in e_{1,1}\mathbb M_d(L(\Gamma_i))p_i$ such that $\omega_i\omega_i^*=e_{1,1}$ and $\mathcal H_i=\mathbb C\omega_i$. In particular, we get that $1=(\uptau\otimes\text{Tr})(e_{1,1})\leqslant (\uptau\otimes\text{Tr})(p_i)=t_i$. Since this holds for all $i\in\overline{1,2}$ and $t_1t_2=1$, we get that $t_1=t_2=1$. Thus, we may assume that $d=1$ and $p_1=p_2=1$.

Let $i\in\overline{1,2}$. Then $\omega_i\in\mathcal U(L(\G_i))$, and since $\mathbb C\omega_i$ is $\alpha_i(\Sigma)$-invariant, we can find a character $\eta_i:\Sigma\ra\mathbb T$ such that $u_{\pi_i(g)}\omega_i\sigma_i(g)^*=\eta_i(g)\omega_i$ and thus $u_{\pi_i(g)}=\eta_i(g)\omega_i\sigma_i(g)\omega_i^*$, for all $g\in\Sigma$. Therefore, if we put $\omega=\omega_1\otimes\omega_2\in L(\G_1)\bar{\otimes}L(\G_2)=M$, then $u_{\pi_i(g)}=\eta_i(g)\omega\sigma_i(g)\omega^*$, for all $g\in\Sigma$.
Denote $u=\omega v^*\in\mathcal U(M)$. Since $L(\La_i)=vL(\G_i)v^*$, we get that $uL(\La_i)u^*=\omega L(\G_i)\omega^*=L(\G_i)$. Moreover, recalling that $v_{\rho_i(\delta(g))}=v\sigma_i(g)v^*$, we get that $u_{\pi_i(g)}=\eta_i(g)uv_{\rho_i(\delta(g))}u^*$, for all $g\in\Sigma$. This implies that $\mathbb T\Sigma=u\mathbb T\Delta u^*$, which finishes the proof.
\hfill$\blacksquare$

Before proving Corollary \ref{B}, we also need the following elementary result:

\begin{lem}\label{super2}
Let $\G$ be an icc group and put $M=L(\G)$. Let  $\Sigma<\G$ be a subgroup such that the centralizer in $\G$ of any finite index subgroup of $\Sigma\cap g\Sigma g^{-1}$ is trivial, for every $g\in\Gamma$. 

If $\Delta<\La$ are countable groups such that $M=L(\La)$ and $\mathbb T\Sigma=\mathbb T\Delta$, then $\mathbb T\G=\mathbb T\La$.
\end{lem}

{\it Proof.} 
In order to prove the lemma, it suffices to show that $\G\subseteq\mathbb T\La$. 

To this end, let $g\in\Gamma$ and put $u:=u_g$. 
Define $\Delta_0:=\Delta\cap \mathbb Tu\Delta u^*$. Then  we have that $\mathbb T\Delta_0=\mathbb T\Delta\cap\mathbb Tu\Delta u^*=\mathbb T\Sigma\cap\mathbb Tg\Sigma g^{-1}=\mathbb T(\Sigma\cap g\Sigma g^{-1})$ and there are a homomorphism $\delta:\Delta_0\rightarrow\Delta$ and a character $\eta:\Delta_0\ra\mathbb T$ such that $u^*v_hu=\eta(h)v_{\delta(h)}$, for all $h\in\Delta_0$. Let $k_1,k_2\in\La$ such that $\uptau(uv_{k_1}^*)\not=0$ and $\uptau(uv_{k_2}^*)\not=0$. Then $\{hk_1\delta(h)^{-1}|h\in\Delta_0\}$ and $\{hk_2\delta(h)^{-1}|h\in\Delta_0\}$ are finite, and hence there is a finite index subgroup $\Delta_1<\Delta_0$ such that $hk_1\delta(h)^{-1}=k_1$ and $hk_2\delta(h)^{-1}=k_2$, for all $h\in\Delta_1$. From this a basic calculation shows that $k:=k_1k_2^{-1}$ commutes with $\Delta_1$. 
Thus, $v_k$ commutes with $\{v_h|h\in\Delta_1\}$ and hence with $\{u_g|g\in\Sigma_1\}$, where $\Sigma_1<\Sigma\cap g\Sigma g^{-1}$ is the finite index subgroup such that $\mathbb T\Sigma_1=\mathbb T\Delta_1$.
The assumption from the hypothesis implies that $v_k\in\mathbb C1$, hence $k=e$ and $k_1=k_2$. Since this holds for every $k_2,k_2\in\La$ in the support of $u$, we conclude that $u\in\mathbb T\La$.
\hfill$\blacksquare$

\subsection*{Proof of Corollary \ref{B}} Let $\Gamma_0$ be an icc, non-amenable, bi-exact group, and $\Sigma_0<\Gamma_0$ be an icc, amenable subgroup. Assume that 
(1) $[\Sigma_0:\Sigma_0\cap g\Sigma_0 g^{-1}]=\infty$, for every $g\in\Gamma_0\setminus\Sigma_0$, and (2) the centralizer in $\G_0$ of any finite index subgroup of $\Sigma_0\cap g\Sigma_0 g^{-1}$  is trivial, for every $g\in\Gamma_0$. Note that (2)  implies that the centralizer of any finite index subgroup of $\Sigma_0$ in $\Gamma_0$ is trivial, or, equivalently, that $\{ghg^{-1}|g\in\Sigma_0\}$ is infinite, for all $h\in\Gamma_0\setminus\{e\}$.

Put $\Gamma_i^j=\Gamma_0$ and $\Gamma_i=\Gamma_i^1\times\Gamma_i^2$, for all $i,j\in\overline{1,2}$. Let $\Sigma=\{(g,g)|g\in\Sigma_0\}<\Gamma_1\cap\Gamma_2$.
Define $\G=\Gamma_1\ast_{\Sigma}\Gamma_2$ and $M=L(\G)$.
Let $\La$ be a countable group such that $M=L(\La)$. 

We first use Theorem \ref{A}.
Let $h=(h_1,h_2)\in\Gamma_i=\Gamma_0\times\Gamma_0$ such that $[\Sigma:\Sigma\cap h\Sigma h^{-1}]<\infty$. Then $h_1h_2^{-1}\in\Gamma_0$ commutes with a finite index subgroup of $\Sigma_0$, and thus by (2) we get that $h_1=h_2$.
Further, it follows that $[\Sigma_0:\Sigma_0\cap h_1\Sigma_0 h_1^{-1}]<\infty$, which by (1) forces $h_1\in\Sigma_0$. Hence $h=(h_1,h_1)\in\Sigma$.
We may thus apply Theorem \ref{A} to deduce that $\Lambda=\Lambda_1\ast_{\Delta}\Lambda_2$ and that, after unitary conjugacy, we have $L(\Sigma)=L(\Delta)$ and $L(\G_i)=L(\La_i)$, for all $i\in\overline{1,2}$.

Since $\{ghg^{-1}|g\in\Sigma_0\}$ is infinite, for all $h\in\Gamma_0\setminus\{e\}$, we are in position to apply Theorem \ref{super}. Thus, we deduce the existence of a decomposition $\La_i=\La_i^1\times\La_i^2$ and a unitary $u_i\in L(\G_i)$ such that $\mathbb T\Sigma=u_i\mathbb T\Delta u_i^*$ and $L(\Gamma_i^j)=u_iL(\La_i^j)u_i^*$, for all $i,j\in\overline{1,2}$. 

In particular, we have that $\mathbb T\Sigma_0=u_i\mathbb T \rho_i^j(\Delta)u_i^*$, where we consider the canonical embedding $\Sigma_0<\Gamma_i^j$ and projection  $\rho_i^j:\Lambda_i\ra\Lambda_i^j$.
By using condition (2) again, Lemma \ref{super2} implies that $\mathbb T\Gamma_i^j=u_i\mathbb T\Lambda_i^ju_i^*$, for all $i,j\in\overline{1,2}$. Thus, we have that $\mathbb T\Gamma_i=u_i\mathbb T\Lambda_iu_i^*$, for all $i\in\overline{1,2}$.

Finally, put $u=u_1u_2^*\in\mathcal U(M)$. Then $\mathbb T\Sigma=u\mathbb T\Sigma u^*$, hence we can find an isomorphism $\delta:\Sigma\ra\Sigma$ and a character $\eta:\Sigma\rightarrow\mathbb T$ such that $u_{\delta(g)}=\eta(g)uu_gu^*$, for all $g\in\Sigma$. Let $k_1,k_2\in \G$ such that $\uptau(uu_{k_1}^*)\not=0$ and $\uptau(uu_{k_2}^*)\not=0$. Then  $\{\delta(g)k_1g^{-1}|g\in\Sigma\}$ and $\{\delta(g)k_2g^{-1}|g\in\Sigma\}$ are finite, hence there is a finite index subgroup $\Sigma_1<\Sigma$ such that $\delta(g)k_1g^{-1}=k_1$ and $\delta(g)k_2g^{-1}=k_2$, for all $g\in\Sigma_1$. 
Since $[\Sigma:\Sigma\cap h\Sigma h^{-1}]=\infty$, for all $h\in\Gamma_i\setminus\Sigma$ and $i\in\overline{1,2}$, we deduce that $k_1,k_2\in\Sigma$.
But then $k:=k_2^{-1}k_1\in \Sigma$ commutes with a finite index subgroup $\Sigma_1<\Sigma$.  Since $\Sigma_0$ is icc, $k=e$, thus $k_1=k_2\in\Sigma$.
Since this holds for any $k_1,k_2\in \G$ in the support of $u$, we derive that $u\in\mathbb T\Sigma$.  

Thus, since $u_1=uu_2$ and  $\mathbb T\Gamma_2=u_2\mathbb T\Lambda_2u_2^*$, we get that $$u_1^*\mathbb T\Gamma_2 u_1=u_2^*(u^*\mathbb T\Gamma_2u)u_2=u_2^*\mathbb T\Gamma_2u_2=\mathbb T\Lambda_2.$$
Since we also have $\mathbb T\Gamma_1=u_1\mathbb T\La_1u_1^*$, we conclude that $\mathbb T\G=u_1\mathbb T\La u_1^*.$ This finishes the proof.
\hfill$\blacksquare$

\section{Proof of Corollary \ref{C}}
Let $\Gamma=\Gamma_1*_{\Sigma}\Gamma_2$ be as in Corollary \ref{B}, where we denote $\Gamma_1=\Gamma_2=\Gamma_0\times\Gamma_0$. Let $\theta:C^*_r(\Gamma)\rightarrow C^*_r(\Lambda)$ be a $*$-isomorphism, for some countable group $\Lambda$. Denote by $\uptau:L(\Lambda)\rightarrow\mathbb C$ the canonical trace and view $C^*_r(\Lambda)\subseteq L(\Lambda)$.
Then $\rho:=\uptau\circ\theta:C^*_r(\Gamma)\rightarrow\mathbb C$ is a tracial state.

We claim that if $g\in\Gamma\setminus\{e\}$, then $\rho(u_g)=0$. To this end, we will show that there exist $a,b\in\Gamma$ such that $b^3\not=e$ and $\{aga^{-1}, b\}$ freely generate a subgroup of $\Gamma$. 

First, assume that $g\in\Sigma\setminus\{e\}$. Then $g=(g_0,g_0)$, for some $g_0\in\Sigma_0\setminus\{e\}$. Since $\Gamma_0$ is icc, we can find $a_0\in\Gamma_0$ such that $a_0$ does not commute with $g_0$. If we put $a=(a_0,e)\in\Gamma_1=\Gamma_0\times\Gamma_0$, then $aga^{-1}=(a_0g_0a_0^{-1},g_0)\in\Gamma_1\setminus\Sigma$. On the other hand, we can find $b_0\in\Gamma_0\setminus\Sigma_0$ such that $b_0^3\not=e$.   Granting this and letting $b=(b_0,e)\in \Gamma_2\setminus\Sigma$, we have that $\{aga^{-1},b\}$ freely generate a subgroup of $\Gamma$.
Now, if we cannot find such a $b_0$, we would have that $b_0^2=e$, for all $b_0\in\Gamma_0\setminus\Sigma_0$.
Thus, if $x,y\in\Gamma_0\setminus\Sigma_0$ are such that $x\Sigma_0\not=y\Sigma_0$, then $x^2=y^2=(x^{-1}y)^2=e$, which implies that $x,y$ commute. Thus, $x\Sigma_0$ and $y\Sigma_0$ commute, which would give that $\Sigma_0$ is abelian, a contradiction.

Secondly, assume that $g\in\Gamma\setminus\Sigma$.  Let $g=g_1g_2...g_k$ be the reduced form on $g$. 
Then the reduced form of $g^n$ begins and ends with $g_1^{\pm}$ or $g_k^{\pm}$, for every $n\in\mathbb Z\setminus\{0\}$.
Let $a\in\Gamma_1\setminus\Sigma$ be such that $a\notin\{g_1^{\pm},g_k^{\pm}\}$. Then the reduced form of $(aga^{-1})^n=ag^na^{-1}$ begins with $a$ and ends with $a^{-1}$, for every $n\in\mathbb Z\setminus\{0\}$. As in the previous paragraph, let $b\in\Gamma_2\setminus\Sigma$ such that $b^3\not=e$. Then it is clear that $\{aga^{-1},b\}$ freely generate a subgroup of $\Gamma$.

Thus, if $\Delta_1,\Delta_2,\Delta<\Gamma$ denote the subgroups respectively generated by $\{aga^{-1}\},\{b\},\{aga^{-1},b\}$, then $\Delta=\Delta_1*\Delta_2$. Since $|\Delta_1|\geqslant 2$ and $|\Delta_2|\geqslant 3$, by Powers' work \cite{Po75} and its  extension \cite{PS79} we get that $C^*_{r}(\Delta)$ has a unique tracial state. Viewing $C^*_r(\Delta)\subseteq C^*_r(\Gamma)$ in the natural way, we conclude that $\rho(u_h)=0$, for all $h\in\Delta\setminus\{e\}$. Thus, $\rho(u_g)=\rho(u_{aga^{-1}})=0$, which proves the claim.
Note that one can alternatively prove the claim by showing that the amenable radical of $\Gamma$ is trivial and applying \cite[Theorem 1.3]{BKKO14}.

Finally, the claim implies that $\rho$ is the restriction of the canonical trace of $L(\Gamma)$ to $C_r^*(\Gamma)$. Thus, $\theta$ is trace preserving and hence it extends to a $*$-isomorphism $\theta:L(\Gamma)\rightarrow L(\Lambda)$. The conclusion now follows from Corollary \ref{B}. 
\hfill$\blacksquare$

\end{document}